\documentclass[notitlepage,twoside,a4paper]{amsart}
\usepackage{mathrsfs}
\usepackage{amsmath,amssymb,enumerate}
\usepackage{epsfig,fancyhdr,color}
\usepackage{epstopdf}
\usepackage{amssymb}
\usepackage{amsmath,amsthm}
\usepackage{latexsym}
\usepackage{amscd}
\usepackage{psfrag}
\usepackage{graphicx}
\usepackage{epsf}
\usepackage[latin1]{inputenc}
\usepackage[all]{xy}
\usepackage{tikz}
\usepackage{prettyref}
\usepackage{subfigure}
\usepackage{float}
\usepackage{color}
\usepackage{array}
\usepackage{hyperref}
\usepackage{cite}
\usepackage[T1]{fontenc}

\usepackage[totalwidth=16cm,totalheight=24cm]{geometry}

\newcommand{\III}{I\hspace{-0.1cm}I\hspace{-0.1cm}I}

\DeclareMathOperator{\area}{Area}

\newtheorem{theorem}{\rm\bf Theorem}[section]
\newtheorem*{theorem*}{\rm\bf Theorem A}
\newtheorem*{theorem*'}{\rm\bf Theorem A*}
\newtheorem*{theorem**}{\rm\bf Theorem B}
\newtheorem*{theorem***}{\rm\bf Theorem C}
\newtheorem*{theorem****}{\rm\bf Theorem D}
\newtheorem*{question*}{\rm\bf Question}
\newtheorem*{remark*}{\rm\bf Remark}

\newtheorem{proposition}[theorem]{\rm\bf Proposition}
\newtheorem{conjecture}[theorem]{\rm\bf Conjecture}
\newtheorem{corollary}[theorem]{\rm\bf Corollary}
\newtheorem{definition}[theorem]{\rm\bf Definition}

\newtheorem{question}[theorem]{\rm\bf Question}

\newtheorem{lemma}[theorem]{\rm\bf Lemma}
\newtheorem{sublemma}[theorem]{\rm\bf Sublemma}
\newtheoremstyle{named}{}{}{\itshape}{}{\bfseries}{.}{.5em}{#1 \thmnote{#3}}
\theoremstyle{named}

\newcommand{\R}{{\mathbb R}}

\newcommand{\CH}{\rm{CH}}

\newcounter{notes}%

\def\interieur#1{\mathord{\mathop{\kern 0pt #1}\limits^\circ}}

\title[]{The induced metric and bending lamination  on the boundary of convex hyperbolic 3-manifolds}

\author{Abderrahim Mesbah}
\address{Abderrahim Mesbah \newline
Universit\'e Paris-Saclay, Laboratoire de Math\'ematiques d'Orsay,\\ 91405 Orsay
, France\\
}
\email{abderrahim.mesbah@universite-paris-saclay.fr}

\date{v0, \today}

\begin{document}

\maketitle

\begin{abstract}
Let $S$ be an oriented closed surface of genus at least two, and let $M = S \times (0,1)$. Suppose that $h$ is a Riemannian metric on $S$ with curvature strictly greater than $-1$, $h^{*}$ is a Riemannian metric on $S$ with curvature strictly less than $1$, and every contractible closed geodesic with respect to $h^{*}$ has length strictly greater than $2\pi$. Let $\mu$ be a measured lamination on $S$ such that every closed leaf has weight strictly less than $\pi$. Then, we prove the existence of a convex hyperbolic metric $g$ on the interior of $M$ that induces the Riemannian metric $h$ (respectively $h^{*}$) as the first (respectively third) fundamental form on $S \times \left\{ 0\right\}$ and induces a pleated surface structure on $S \times \left\{ 1\right\}$ with bending lamination $\mu$. This statement remains valid even in limiting cases where the curvature of $h$ is constant and equal to $-1$. Additionally, when considering a conformal class $c$ on $S$, we show that there exists a convex hyperbolic metric $g$ on the interior of $M$ that induces $c$ on $S \times \left\{ 0\right\}$, which is viewed as one component of the ideal boundary at infinity of $(M,g)$, and induces a pleated surface structure on $S \times \left\{ 1\right\}$ with bending lamination $\mu$. Our proof differs from previous work by Lecuire for these two last cases. Moreover, when we consider a lamination which is small enough, in a sense that we will define, and a hyperbolic metric, we show that the metric on the interior of M that realizes these data is unique.
\end{abstract}
\section*{introduction}
Let $S$ be a closed surface of genus $g \geq 2$. We denote by $M$ the three-dimensional manifold $S \times \left( 0,1 \right)$, and we denote by $\overline{M}$ the three-dimensional manifold with boundary $S \times \left[ 0,1 \right]$. We also denote by $\mathcal{T}(S)$  the Teichm\"{u}ller space of $S$, which is the space of isotopy classes of hyperbolic structures on $S$, or equivalently, due to uniformization, the space of isotopy classes of conformal structures on $S$.\\
We denote the isometry group of the hyperbolic space $\mathbb{H}^3$ as $Isom(\mathbb{H}^3)$, and the subgroup of isometries that preserve orientation as $Isom^{+}(\mathbb{H}^3)$.\\
Let $\Gamma$ be a subgroup of $Isom^{+}(\mathbb{H}^3)$. We say that $\Gamma$ is Kleinian if it acts properly discontinuously on the hyperbolic space.\\
\subsection{Convex cocompact hyperbolic manifolds}
\begin{figure}
    \centering
    \includegraphics[scale=0.6]{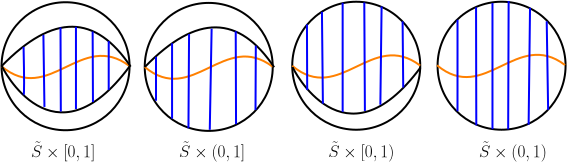}
\caption{The shadowed region is denoted by $\Omega$. The number of boundary components of the quotient depends on the number of boundary components of $\Omega$. We always have two boundary components when the ideal boundary is included. In the case where the interval $(0,1)$ is open on both sides, $\Omega$ is equal to $\mathbb{H}^3$.}
    \label{depend-fig}
\end{figure}

We say that a hyperbolic manifold (possibly with boundary) is geodesically convex if every arc in the manifold is homotopic, relative to its endpoints and within the manifold, to a geodesic arc contained entirely in the manifold. Through this paper, convex hyperbolic manifolds are always assumed to be complete (unless mentioned otherwise).\\
Thurston has shown (see \cite[Proposition 8.3.2]{thurston2022geometry}) that if $N$ is a geodesically convex hyperbolic manifold, then there exists a closed convex subset $\Omega$ of $\mathbb{H}^3$ and a Kleinian group $\Gamma$ such that $N$ is isometric to $\Omega / \Gamma$. Moreover, if the quotient of the convex hull of the limit set of $\Gamma$ by $\Gamma$ is compact, we say that $N$ is a hyperbolic convex cocompact manifold (see Section 2 for a definition of the limit set and convex hull).\\
We will focus on convex cocompact hyperbolic manifolds of the form $\Omega / \Gamma$, where $\Gamma$ is a Kleinian group that has limit set equal to a Jordan curve and $\Omega$ is a closed convex subset of $\mathbb{H}^3$ that is invariant under $\Gamma$, assume moreover that $\Gamma$ is isomorphic to  $\pi_{1}(S)$, the fundamental group of $S$. Then, according to \cite[Proposition 8.7.2]{thurston2022geometry}, the manifold $\Omega / \Gamma$ is diffeomorphic to either $S \times [0,1]$, $S \times (0,1]$, $S \times [0,1)$, or $S \times (0,1)$, depending on the number of boundary components $\Omega$ has. Since convex hyperbolic manifolds are considered to be complete here, unless mentioned otherwise, we interpret an open interval as having an ideal boundary, and a closed interval as having a boundary component of $\Omega$ inside $\mathbb{H}^3$ (see Figure \ref{depend-fig}). Equivalently, if $g$ is a convex cocompact hyperbolic metric on one of these manifolds, then the resulting manifold is identified with $\Omega / \Gamma$, where $\Gamma$ is a Kleinian group with a limit set equal to a quasi-circle, and $\Omega$ is a closed convex subset of $\mathbb{H}^3$ invariant under the action of $\Gamma$.\\

If $\Omega$ is the entire $\mathbb{H}^3$, we refer to $\mathbb{H}^3 / \Gamma$ as a quasi-Fuchsian (see Section 2 for more details). Note that in all cases, a hyperbolic convex cocompact metric on one of the four manifolds mentioned above is embedded in a unique quasi-Fuchsian manifold, which arises from the fact that $\Omega / \Gamma$ is embedded in $\mathbb{H}^{3} / \Gamma$.\\
\subsection{The boundary of convex cocompact hyperbolic manifolds}
If the manifold is complete and the boundary is open on one side or both sides, such as $S \times [0,1)$, $S \times (0,1)$, $S \times [0,1)$, the convex, cocompact, complete hyperbolic manifold can still be embedded in $S \times [0,1]$ by adding the ideal boundary of $\Omega$ minus the limit set of $\Gamma$ (see Section 2 for details). But when we extend in this way, the metric cannot be extended to this ideal boundary (because the hyperbolic metric cannot be extended to the ideal boundary of $\mathbb{H}^{3}$), but it induces a conformal structure (see Section 2 for details). In this case, we refer to the boundary as the boundary at infinity. It can be $S \times \left\{0\right\}$, $S \times \left\{1\right\}$ or $S \times \left\{0,1\right\}$ (see Section 2 for more details). In all cases, the convex, cocompact hyperbolic metric on $S \times (0,1)$ induces data on the boundary. It can be a conformal structure when we take the ideal boundary as explained above, or it can be an induced metric which may not be smooth, as we will see, but in all cases this metric has a Gaussian curvature bigger than or equal to $-1$. If a boundary component of the hyperbolic convex cocompact manifold $S \times [0,1]$ is not ideal and is smoothly embedded in a quasi-Fuchsian manifold, then the induced metric has a Gaussian curvature strictly greater than $-1$. In this case, we refer to the induced metric as the first fundamental form. Additionally, we can define another metric on this boundary that measures how much the embedded surface is curved inside the quasi-Fuchsian manifold. This is referred to as the third fundamental form (see Section 1.4 for details on the third fundamental form). If both boundary components of the hyperbolic convex cocompact manifold $S \times [0,1]$ are smoothly embedded in a quasi-Fuchsian manifold, we say that the hyperbolic convex cocompact manifold $S \times [0,1]$ is strictly convex. In other words, the hyperbolic convex cocompact manifold $\overline{M}$ is considered strictly convex when there are no geodesic arcs that are fully contained within its boundary.\\
\subsection{Convex core of quasi-Fuchsian manifolds}
Let $\Gamma$ be a Kleinian group, $\Lambda_{\Gamma}$ be its limit set, and $CH(\Lambda_{\Gamma})$ be the convex hull of the limit set (inside $\mathbb{H}^{3}$). If the limit set of a Kleinian group is a quasi-circle, which is not a geometric circle, then the quotient $CH(\Lambda_{\Gamma}) / \Gamma$ is a compact manifold homeomorphic to $S \times \left[0,1\right]$ (and a totally geodesic surface homeomorphic to $S$ if $\Lambda_{\Gamma}$ is a geometric circle), endowed with a convex cocompact hyperbolic metric. Denoting $\mathbb{H}^{3} / \Gamma$ as $Q$ (a quasi-Fuchsian manifold), we refer to $CH(\Lambda_{\Gamma}) / \Gamma$ as $C(Q)$, calling it the convex core of $Q$. $C(Q)$ is the smallest complete geodesically convex (it contains all the simple closed geodesics) submanifold of $Q$, as detailed in \cite[Proposition 8.1.2]{thurston2022geometry}. We see $C(Q)$ as a submanifold of $Q$. The boundary of $C(Q)$ consists of two components that are  pleated surfaces (see Section 1.3 for the definition). The hyperbolic metric of $Q$ induces a path metric on the boundary of $C(Q)$ (even if the boundary is not smoothly embedded. see \cite[Section 8.5]{thurston2022geometry} and \cite[Chapter II.1]{canary2006fundamentals}), which is isometric to a hyperbolic metric, and the pleating induces a measured lamination referred to as the bending lamination (see Section 1.2 for the definition of measured lamination and Section 2.2 for an explanation of how pleating induces a measured lamination).\\
\subsection{Various types of boundary data}
As we can see, a hyperbolic convex cocompact metric on $S \times (0,1)$ can induce several types of data on the boundary (that is, $S \times \left\{0,1\right\}$). In particular, this data can include a metric with curvature strictly greater than $-1$, which is referred to as the first fundamental form or simply the induced metric. It can also include a metric with curvature strictly smaller than $1$, known as the third fundamental form and usually denoted as $\III$. Furthermore, in the case where the boundary is ideal (meaning it lifts to the ideal boundary of $\mathbb{H}^3$), this data is a conformal structure. Alternatively, when the boundary is not smoothly embedded but is a pleated surface (namely, the boundary of the convex core), the data can be a hyperbolic metric, or a measured lamination known as the bending lamination.\\
The study of convex cocompact hyperbolic metrics on  $M$ , and more precisely the correspondences between the metric $g$ on $M$ and the boundary data on $M$ , has been the subject of significant research interest.\\
Let's consider a quasi-Fuchsian manifold $Q$, recall that this manifold has two ideal boundary components, and then it induces a conformal structure on $S \times \left\{0,1\right\}$. It was shown by Bers \cite{bams/1183523461} (see Theorem \ref{bers}) that there is a homeomorphism (actually a biholomorphism if we consider the complex structure induced by the character variety) between the space of conformal structures on $\partial M$ and the space of quasi-Fuchsian structures on $M$.\\
\subsection{Foliation by $k$-surfaces}
A key result about quasi-Fuchsian manifolds is that for any quasi-Fuchsian manifold $Q$, each connected component of $Q\setminus C(Q)$ is foliated by $k$-surfaces (see \cite{labourie1991probleme}). That is, $Q\setminus C(Q) = \bigcup\limits_{k\in(-1,0)} S^{\pm}_k$, where the signs $\pm$ on $S^{+}_k$ and $S^{-}_k$ indicate to which connected component the surface belongs. The induced metric on each of $S^{+}_k$ and $S^{-}_k$ has a constant Gaussian curvature equal to $k$ (see Section 2.3).\\
\subsection{Parametrization of convex cocompact hyperbolic manifolds by boundary data}
Labourie, in \cite{labourie1992metriques}, showed that for any metric $h$ on $S \times \left\{0,1\right\}$ with Gaussian curvature strictly greater than $-1$, one can find a unique hyperbolic convex cocompact metric on $M$ that induces $h$ on the boundary. As previously explained, $(\overline{M},g)$ is embedded in the interior of a quasi-Fuchsian manifold. Later, Schlenker, in \cite{MR2208419}, established the uniqueness  of the hyperbolic convex cocompact metric $g$ on $M$ that induces $h$ on the boundary of $\overline{M}$. He also showed a similar statement regarding the uniqueness and existence of the third fundamental form. Specifically, any metric $h^{*}$ on $S \times \left\{0,1\right\}$ with Gaussian curvature strictly smaller than $1$, and in which every closed geodesic has length strictly greater than $2\pi$, can be induced as the third fundamental form on the boundary of $\overline{M}$ by a unique hyperbolic convex cocompact metric $g$ on $M$, where $(\overline{M},g)$ is once again embedded in the interior of a quasi-Fuchsian manifold. If we assume that $h$ or $h^{*}$ has constant Gaussian curvature, then $S \times \left\{0,1\right\}$ corresponds to some $k$ surfaces, while we observe $(\overline{M},g)$ is embedded within a quasi-Fuchsian manifold.\\   Moreover, in \cite[Theorem 2.6]{MR2208419} the author shows that for each $k\in(-1,0)$, the map that maps a quasi-Fuchsian structure on $M$ to the induced metric on $S_k^- \bigsqcup S_k^+$ is a homeomorphism (in fact, a diffeomorphism. See \cite{MR2208419} for details). For any quasi-Fuchsian manifold $Q$, when $k$ goes to $-1$, the induced metrics on $S^{\pm}_k$ converge to the induced metrics on $\partial^{\pm}C(Q)$ (the boundary of the convex core) in the length spectrum, and the third fundamental forms induced on $S^{\pm}_k$ converge in the length spectrum to the bending lamination of $\partial^{\pm}C(Q)$ (see Theorem 2.7). Thurston has conjectured that the map that maps a quasi-Fuchsian structure on $M$ to the induced metric on the boundary of its convex core is a homeomorphism. Even though it is known that any two hyperbolic metrics on $S$ can be realized as the induced metrics on the boundary of the convex core of some quasi-Fuchsian manifold $Q$ (up to isotopy), the conjecture remains unproven, as injectivity has not been established, except in certain special cases. For example, Prosanov proved the uniqueness of the manifold for an open and dense subset of the deformation space of quasi-Fuchsian manifolds (see \cite{Prosanov}). Thurston made a similar conjecture for bending laminations, namely that the map that associates to a quasi-Fuchsian structure on $M$  the bending lamination on the boundary of its convex core is a homeomorphism. In \cite{bonahon2004laminations}, the authors characterize the measured laminations on $\partial M$ that can be realized (up to isotopy) as the bending lamination of the convex core. Later, Dular and Schlenker (see \cite{Dular-schlenker}) showed that if a pair of measured laminations on $S$ arises as the bending laminations of the boundary of the convex core of a quasi-Fuchsian manifold $Q$, then $Q$ is uniquely determined.\\
 In the paper \cite{chen2022geometric}, the authors introduced a new parametrization of quasi-Fuchsian structures on $M$. This parametrization is obtained by considering the induced third fundamental form on $S^+_k$ and the induced metric on $S^-_k$, or by taking the third fundamental form on $S^+_k$ and the conformal structure on $S\times  \left\{ 0\right\}$ (see Theorem 2.8). The main result of this paper is the extension of the surjectivity statements to the boundary of the convex core.\\
 In this paper, we consider only convex hyperbolic metrics on $M$ that are cocompact. In particular, a convex hyperbolic metric on $M$ is a metric $g$ such that $(M, g)$ is isometrically embedded in a quasi-Fuchsian manifold (which must be unique), and the image is geodesically convex. For more details, we refer to \cite[Section 8.3]{thurston2022geometry}.
\subsection{Main theorems of the paper}
 
 Our main theorems are as follows.

\begin{theorem*}
Let $h$  be Riemannian metric on $S$, and denote its curvature by $k_h$ (which is not necessarily constant). We assume that $-1 < k_{h}$. Let $\mu$ be a measured lamination on $S$ such that every closed leaf has weight strictly smaller than $\pi$. Then, there exists a convex hyperbolic metric $g$ on $M = S \times (0,1)$, the interior of $\overline{M} = S \times \left[0,1\right]$, such that:
\begin{itemize}
    \item $g$ induces a metric on $S\times \left\{ 0 \right\}$ which is isotopic to $h$.
    \item $g$ induces on $S \times \left\{ 1\right\}$ a pleated surface structure in which its bending lamination is $\mu$.
\end{itemize}
Furthermore, when we see $\overline{M}$ as an embedded manifold in a quasi-Fuchsian manifold $Q$, then the boundary $S\times \left\{ 0 \right\}$ is smoothly embedded.
\end{theorem*}

We prove a similar statement for third fundamental forms.

\begin{theorem*'}
Let $h^{*}$  be a Riemannian metric on $S$, and denote its curvature by $k_{h^*}$ (which is not necessarily constant). We assume that $k_{h^{*}} < 1$. Moreover, we assume that every contractible closed geodesic with respect to $h^{*}$ has length strictly bigger than $2\pi$. Let $\mu$ be a measured lamination on $S$ such that every closed leaf has weight strictly smaller than $\pi$. Then, there exists a convex hyperbolic metric $g$ on $M=S\times(0,1)$, the interior of $\overline{M}=S\times[0,1]$, such that:
\begin{itemize}
    \item $g$ induces a third fundamental form on $S\times \left\{ 0\right\}$ which is isotopic to $h^{*}$.
    \item $g$ induces on $S \times \left\{ 1\right\}$ a pleated surface structure in which its bending lamination is $\mu$.
\end{itemize}
Furthermore When we see $\overline{M}$ as an embedded manifold in a quasi-Fuchsian manifold $Q$, then the boundary $S\times \left\{ 0 \right\}$ is smoothly embedded.
\end{theorem*'}

However, it is currently unknown whether the convex hyperbolic metric $g$ on $M$ is unique under the
conditions stated in Theorem A.
\begin{question}
Is the metric $g$ given in Theorem A or in Theorem A* unique (up to isotopy) under the hypothesis of the theorems ?
\end{question}

\begin{figure}
    \includegraphics[scale=0.6]{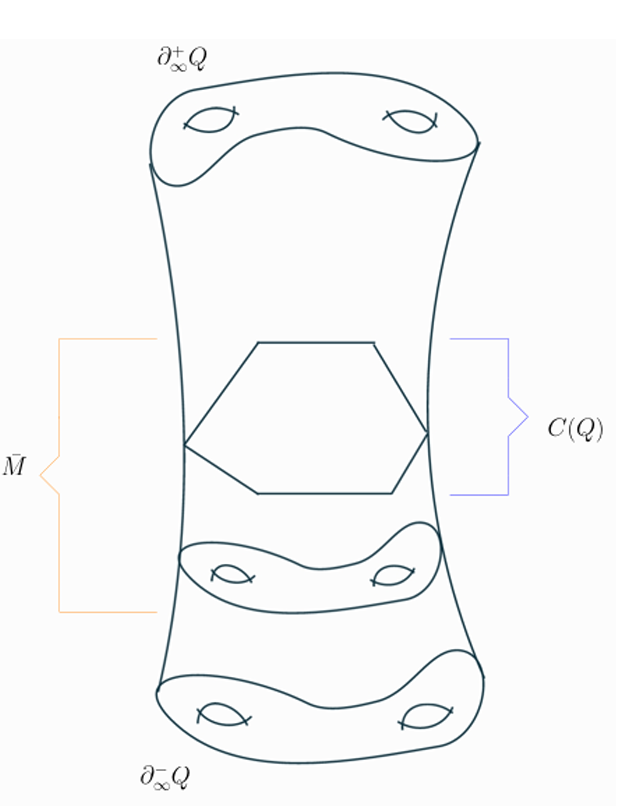}
\caption{The manifold $\overline{M}$ is embedded in a quasi-Fuchsian manifold $Q$. One of its boundaries is smoothly embedded in $Q$, while the other one is a boundary component of $C(Q)$.}
    \label{The_theorem}
\end{figure}

For a quasi-Fuchsian manifold $Q$, we will denote by $\partial^{+}C(Q)$ (resp. $\partial^{-}C(Q)$) the upper (resp. lower) boundary component of the convex core of $Q$, and we will denote by $\partial^{+}_{\infty}Q$ (resp. $\partial^{-}_{\infty}Q$) the upper (resp. lower) boundary component at infinity of $Q$ (see Figure \ref{The_theorem}).\\
The proof leads us to show that the map that associates to a quasi-Fuchsian manifold the induced metric on $\partial^- C(Q)$ and the bending lamination on $\partial^+ C(Q)$ (call it $Mix_{-1}$), and the map that associates to a quasi-Fuchsian manifold $Q$ the induced conformal structure on $\partial^{-}_{\infty}Q$ and the bending lamination on $\partial^+ C(Q)$ (call it $Mix_{\infty, -1}$) are proper (see Section 2.4 for more details).

So, another theorem that we show in this paper is the following:
\begin{theorem**}\label{hope}
 The maps $Mix_{-1}$ and $Mix_{\infty,-1}$ are proper.   
\end{theorem**}
In other words, Theorem B states the following. Consider a sequence $(Q_n)_{n \in \mathbb{N}}$ of quasi-Fuchsian manifolds. If the sequence $(h_n)_{n \in \mathbb{N}}$ of the hyperbolic metrics on $\partial^{-} C(Q_n)$ (up to isotopy) converges to a hyperbolic metric $h$, and if $(B^{+}_{n})_{n \in \mathbb{N}}$, the sequence of bending laminations on $\partial^+ C(Q_n)$, converges to $B_{\infty}$, where $B_{\infty}$ is a measured lamination in which each closed leaf of it has weight strictly less than $\pi$, then $(Q_n)_{n \in \mathbb{N}}$ converges, up to extracting a subsequence to a quasi-Fuchsian manifold.\\
The proof of Theorem B employs techniques similar to those used in \cite[Section 3]{lecuire2014convex}, and the proof of this theorem will take the most important part of the paper. The proof consists of showing that the two induced metrics on $\partial^\pm C(Q_n)$ (up to isotopy) are in a compact subset of Teichm\"{u}ller space. Then, we use Theorem \ref{sullivan} to deduce that the conformal structures at the boundary at infinity are in a compact subset of Teichm\"{u}ller space. Finally, we conclude by the theorem of Bers (Theorem 2.1). After having proven Theorem B, we use the parametrization introduced in \cite{chen2022geometric} (see Theorem 2.8) to approximate the metric and the measured lamination that we want to realize by $k$-surfaces. We use Theorem B to show that the sequence of quasi-Fuchsian manifolds that we constructed has a subsequence that converges, and finally, we conclude by Theorem 2.7.\\
As a limiting case, we recover Theorem C, which is stated in \cite{lecuire2006mixing} (in unpublished notes). However, the proof we provide is completely independent from the one given in \cite{lecuire2006mixing}. We believe that our arguments can also be applied to quasi-Fuchsian manifolds with particles. These manifolds were introduced by Thurston (see, for example, \cite{cone}), they share the same geometric structure as quasi-Fuchsian manifolds, except along a finite number of infinite lines where they exhibit cone singularities. Quasi-Fuchsian manifolds with particles have been the subject of recent research (see, for example, \cite{lecuire2014convex}, \cite{moroianu2009quasi}, and \cite{chen-JM-particles}), particularly due to their connection with the study of hyperbolic surfaces with conic singularities. In the convex core of such a manifold, the bending locus consistently avoids the singularities, depending on the conical angles. As a result, there is strong evidence that the arguments used in our proof of Theorem B can be extended to this setting (for more details, see \cite{moroianu2009quasi} and \cite{lecuire2014convex}). This gives us hope that Theorem C also holds for quasi-Fuchsian manifolds with particles. Nevertheless, it remains unknown whether the methods used in Section 4 apply to this case.
  
\begin{theorem***}
 Let $\mu$ be a measured lamination in which every closed leaf has weight strictly smaller than $\pi$, and let $h\in \mathcal{T}(S)$. Then: 
 \begin{itemize}
     \item there exists a quasi-Fuchsian manifold $Q$ such that the bending lamination of $\partial^+ C(Q)$ is $\mu$, and the induced metric on $\partial^- C(Q)$ is isotopic to $h$.
     \item There exists a quasi-Fuchsian manifold $Q$ such that the bending lamination of $\partial^+ C(Q)$ is $\mu$, and the induced conformal structure on the lower boundary at infinity $\partial_{\infty}^{-}Q$ of $Q$ is isotopic to $h$.
 \end{itemize}
\end{theorem***}

A natural question that arises from Theorem C is the uniqueness of the quasi-Fuchsian  manifold that realizes $\mu$ and $h$. We reformulate these questions in the following way.

\begin{question}
  Let $\mu$ be a measured lamination in which every closed curve has weight strictly smaller than $\pi$, and let $h\in T(S)$. Is there a unique quasi-Fuchsian manifold $Q$ that induces $\mu$ as the bending lamination on $\partial^+ C(Q)$ and induces $h$ as the hyperbolic metric on $\partial^- C(Q)$ ?.  
\end{question} 

In Section 6, we will give a partial answer to the last question.

\begin{theorem****}
  For any measured lamination $\mu$ and a hyperbolic metric $h$ on $S$, there exists an $\delta_{h,\mu} > 0$ such that for any $0 < t < \delta_{h,\mu}$, there exists a unique quasi-Fuchsian manifold $Q$, such that the hyperbolic metric on $\partial^- C(Q)$ is $h$ (up to isotopy) and the bending lamination on $\partial^+ C(Q)$ is $t\mu$.  
\end{theorem****}

Even if we don't give an answer on the following question, we believe that a partial answer similar to Theorem D, can be treated by similar arguments as the ones in Section 6.
\begin{question}
  Let $\mu$ be a measured lamination in which every closed curve has weight strictly smaller than $\pi$, and let $h\in T(S)$. Is there a unique quasi-Fuchsian manifold $Q$ that induces $\mu$ as the bending lamination on $\partial^+ C(Q)$ and induces $h$ as the conformal structure on $\partial^{-}_{\infty}Q$? 
\end{question} 
For the sake of clarity, we begin by proving Theorem C in Section 4 before moving on to a proof of Theorem A in Section 5. While the techniques used in both proofs are largely the same, it is especially important to be precise when proving Theorem A*.\\
We will need to introduce some more background in Section 6, we prefer to add it in Section 6 rather than Section 1 for the reader's comfort.

\subsection*{Outline of the paper}
in the first section we recall the necessary background about Teichm\"{u}ller space, measured laminations, pleated surfaces, third fundamental forms, and the duality between the hyperbolic space and de Sitter space.\\
In Section 2, we introduce the quasi-Fuchsian manifolds, and we give some key theorems about them. Namely, we will give the statements about the parametrization by the conformal structure on the boundary at infinity, the parametrization by the induced metrics and the induced third fundamental forms on the k-surfaces, Thurston conjectures about the induced metric and the bending lamination on the boundary of the convex core. We finish the section by giving the parametrization introduced in \cite{chen2022geometric}, and introducing the main theorems that we will prove in this paper.\\
In Section 3, we show that the maps $Mix_{-1}$ and $Mix_{\infty,-1}$ introduced in Section 2.4 are proper. We call this section the "closing lemma," and it will be the crucial part of this paper.\\
In Section 4 we conclude Theorem C by using an approximation by k-surfaces and the parametrization introduced in \cite{chen2022geometric}.\\
In Section 5 we give a proof for the main theorems by using similar techniques as in Section 4.\\
In Section 6 we introduce more background using the differential structure of Teichm\"{u}ller space and quasi-Fuchsian manifolds and the tangentiable structure on the measured laminations to establish Theorem D by using a similar approach to the one used in \cite{bonahon2005kleinian}.

\subsection*{Acknowledgments}

I would like to express my gratitude to my supervisor, Jean-Marc Schlenker, for his invaluable help and unwavering support. I would also like to extend my thanks to Fran\c{c}ois Fillastre for bringing to our attention the fact that Theorem C has already been proven by Lecuire.\\
Likewise, I wish to extend my appreciation to my colleagues from the Department of Luxembourg, particularly Tommaso Cremaschi, for the insightful discussions we had regarding Section 6. Additionally, I am thankful to Mehdi Belraouti and Lamine Messaci for their contributions to the discussions about the background of Section 6.
This work has been supported by the Luxembourg National Research Fund PRIDE/17/1224660/GPS.
\section{PRELIMINARIES}
\subsection{Teichm\"{u}ller space}

For more details see \cite[Section 7]{fathi2021thurston}. Through the paper $S$ will be considered to be a smooth closed connected surface of genus $g \geq 2$. The Teichm\"{u}ller space of $S$, denoted by $\mathcal{T}(S)$, is the set of hyperbolic metrics on $S$ modulo isotopy. That is, two metrics $h_1$ and $h_2$ are identified if and only if there is an isometry $I:(S,h_1) \to (S,h_2)$ which is isotopic to the identity. Thanks to the uniformisation theorem, Teichm\"{u}ller space can also be seen as the space of conformal structures on $S$ modulo isotopy. That is, two conformal structures $c_1$ and $c_2$ are identified if and only if there is a conformal diffeomorphism $f:(S,c_1) \to (S,c_2)$ which is isotopic to the identity.\\
Let $\mathcal{S}$ be the set of free homotopy classes of simple closed curves not homotopic to one point on $S$. There is an embedding of $\mathcal{T}(S)$ into $\mathbb{R}^{\mathcal{S}}_{\geq0}$ by the map, $h \mapsto \Phi(h)$, where $h$ is a hyperbolic metric and $\Phi_{h}$ is the map 
\begin{align*}
 \Phi(h): \mathcal{S} & \to \mathbb{R}_{\geq0}\\
 \alpha & \mapsto l_{h}(\alpha')
\end{align*}
where $\ell_{h}(.)$ is the length function with respect to $h$, and $\alpha'$ is the unique simple closed geodesic in the free homotopy class $\alpha$.\\
This last function $\Phi$ can be defined for any metric with strictly negative curvature, and it is well defined modulo isotopy. That is, if there is an isometry which is isotopic to the identity between $(S,h)$ and $(S,h')$ then $\Phi(h) \equiv  \Phi(h')$. We denote $\Phi(h)$ by $\ell_{h}$.\\ 
Through this paper when we mention a metric of negative curvature we mean its isotopy class, that is, we identify any two metrics that are isometric to each other  via an isometry isotopic to the identity. We say that a sequence $(h_{n})_{n \in \mathbb{N}}$ of metrics of negative curvature converges to a metric $h$ (also with negative curvature) in the length spectrum, if $\ell_{h_{n}}$ converge to $\ell_{h}$ point-wisely.\\ 
The embedding $\Phi$ determines a natural topology on $\mathcal{T}(S)$ in which a sequence of points $h_n$ in $\mathcal{T}(S)$ converge to a point $h$ in $\mathcal{T}(S)$ if and only if $h_n$ converge to $h$ in the length spectrum. (To be more precise, the embedding $\Phi$ is proper, when we endow $\mathcal{T}(S)$ with the quotient topology coming from seeing the hyperbolic metrics on $S$ as tensors, see reference above for more details).
\subsection{Measured laminations}
For more details see \cite[Chapter I.4]{canary2006fundamentals}, \cite{lyubich2001laminations}, \cite[Chapter 8.6]{thurston2022geometry}. Provide $S$ with a hyperbolic metric $h$. We say that $\mu$ is a geodesic lamination on $S$ if it is a closed set that is a disjoint union of simple geodesics. The connected components of $\mu$ are called its leaves, and the connected components of $S \setminus \mu$, together with the leaves of $\mu$, are called strata. Even though the notion of a geodesic lamination seems to depend on the hyperbolic metric given on $S$, it actually depends only on the topology of $S$. Indeed, if we equip $S$ with a new hyperbolic metric $h'$, each leaf $l$ of $\mu$ becomes a quasi-geodesic in $(S,h')$. Then, the lift $\tilde{l}$ of $l$ in $\mathbb{H}^2$ determines a unique geodesic (the geodesic that connects the endpoints of that leaf). Consequently, each leaf $l$ determines a unique geodesic on $(S,h')$, which is the projection of the geodesic that connects the endpoints of $\tilde{l}$.   \\

A transverse measure $\lambda$ on $\mu$ is the assignment of a Radon measure $\lambda_{\kappa}$ on each arc $\kappa$ which is transversal to $\mu$ such that:

\begin{itemize}
    \item If $\kappa'$ is a sub-arc of $\kappa$ then $\lambda_{\kappa}|_{\kappa'}=\lambda_{\kappa'}$.
    \item If $\kappa$ and $\kappa'$ are two arcs transversal to $\mu$, and homotopic relative to $\mu$, then $\lambda_{\kappa}$ and $\lambda_{\kappa'}$ are compatible. More precisely if $H:\left [0,1  \right ]\times \left [0,1  \right ] \to S$, is a homotopy between $\kappa$ and $\kappa'$ relative to $\mu$ (that is, if $H(x,0)$ belongs to $\mu$, then $H(x,t)$ belongs to $\mu$ for any $t$ in $\left [0,1  \right ]$ ), we denote $H_{t} = H(.,t)$ and we assume that $H_{t}$ is a continuous embedding for any $t$. We denote $tr_{\kappa',\kappa} := H_{0}\circ H^{-1}_{1}$ which a homeomorphisme between $\kappa'$ and $\kappa$, then  $\lambda_{\kappa'} = tr_{\kappa',\kappa}^{*}\lambda_{k}$.  
\end{itemize}
Thought this paper, unless there is a confusion, we will refer to both of  $(\mu,\lambda)$ and $\mu$ by $\lambda$. We denote the set of measured laminations on $S$ by $\mathcal{ML}(S)$.\\
Note that we can assign to each transverse arc $\kappa$ a positive number $i(\kappa,\lambda) := \int_{\kappa}d\lambda_{\kappa}$, in fact $i(.,\lambda)$ determines completely the measured lamination $\lambda$, see \cite[Section 6]{BONAHON1997103}.\\
The simplest example of a measured lamination, is in the case when the support of $\lambda$ consists of finitely many closed leaves, and we associate to each leaf an atomic unit mass. In this case we say that $\lambda$ is rational (or discrete). \\ 
Exactly as for Teichm\"{u}ller space, there is a proper embedding of $\mathcal{ML}(S)$ into $\mathbb{R}^{\mathcal{S}}_{\geq0}$ by the map $\lambda \in \mathcal{ML}(S) \mapsto \mathcal{I}(\lambda)$, such that $\mathcal{I}(\lambda)$ is the map
\begin{align*}
\mathcal{I}(\lambda): \mathcal{S} &\to \mathbb{R}_{>0}\\
\alpha &\mapsto i(\alpha',\lambda)
\end{align*}
where $\alpha'$ is the unique closed geodesic in the free homotopy class $\alpha$ (the map $\mathcal{I}$ does not depend on the chosen metric). If $\alpha'$ is not transversal to $\lambda$ then either it belongs to its leaves, or it is disjoint from $\lambda$, in this case we put $\mathcal{I}(\lambda)(\alpha) = 0$. This embedding gives a natural topology on the space of measured laminations in which $(\lambda_{n})_{n \in \mathbb{N}}$ converge to $\lambda$ if and only if  $(i(\cdot,\lambda_{n}))_{n \in \mathbb{N}}$ converge to $i(\cdot,\lambda)$ point-wisely. We call this topology the weak*-topology, for more details about this we refer to \cite[Section 6]{fathi2021thurston}.\\
We say that a sequence of metrics $(h_{n})_{n \in \mathbb{N}}$ of negative curvature converge to a measured lamination $\lambda$ in the length spectrum if for any $\alpha \in \mathcal{S}$, $\ell_{h_{n}}(\alpha)$ converge to $i(\alpha,\lambda)$.

\subsection{Pleated surfaces}
Let $M$ be a hyperbolic 3-manifold. A pleated surface of $M$ is a couple $(S,f)$, where $S$ is a complete hyperbolic surface, and $f:S \to M$ is a map that satisfies the following property: Every point in $S$ is contained in the interior of a geodesic arc of $S$ which is mapped by $f$ to a geodesic arc of $M$.\\
The pleating locus of a pleated surface $(S,f)$ is the set of points of $S$ that are contained in exactly one geodesic arc which is mapped to a geodesic arc of $M$. Since $S$ is connected, the pleated locus of $S$, if non-empty, is a geodesic lamination, see \cite[Section I.5]{canary2006fundamentals}. 
\subsection{Third fundamental form}
Let $M$ be a hyperbolic 3-manifold, and let $S$ be an immersed surface. The Riemannian metric of $M$ when restricted to the tangent bundle of $S$ gives a Riemannian metric on $S$ which is called the induced metric or the first fundamental form, and is denoted by $I$. Let $N$ be a unit normal vector field on $S$, and let $\nabla$ be the Levi-Civita connection of $M$, then the shape operator $B: TS \to TS$ is defined to be $Bu= -\nabla_{u}N$.\\
The third fundamental form $\III$ of $S$ is defined by:
$$\forall x \in S, \ \forall u,v \in T_{x}S, \ \III(u,v) = I(Bu,Bv). $$
The extrinsic curvature $K_{ext}$ of $S$ is the determinant of the shape operator $B$. This quantity is related to $K$, the Gaussian intrinsic  curvature (or sectional curvature) of $S$,  by the equation $K = K_{ext} - 1$.\\
 Using \cite[Proposition 3.12]{MR2328927}, we get that the Gaussian curvature $K^*$ of $S$ endowed with the metric $\III$ is given by:
$$K^* = \frac{K}{K+1}. $$
\subsection{Duality between hyperbolic and de Sitter geometry}
In order to geometrically interpret the third fundamental form $\III$, we need to introduce the de Sitter space $dS^3$ and clarify the connection between $\mathbb{H}^3$ and $dS^3$. Let $\mathbb{R}^{3,1}$ be the vector space $\R^{4}$ endowed with Lorentzian scalar product $\left<\cdot,\cdot \right>_{3,1}$ of signature $(3,1)$.\\ The hyperbolic space can be viewed as $\mathbb{H}^3 := \left\{x = (x_1,x_2,x_3,x_4) \in \mathbb{R}^4 | \left<x,x \right>_{3,1} = -1, x_4 >0  \right\} $ with the metric induced from $\left<\cdot,\cdot \right>_{3,1}$. De Sitter space is defined to be the Lorentzian analogue of the sphere in the euclidean space, that is $dS^3 := \left\{y \in \mathbb{R}^4 | \left<y,y \right>_{3,1} = 1 \right\} $.\\
Let $P \subset \mathbb{H}^3$ be a totally geodesic plane, and let $n_{P}$ be the unit normal vector on $P$. Note that $\left<n_{P},n_{P} \right>_{3,1} = 1$. If we see $n_{P}$ as a point in $\mathbb{R}^{3,1}$, it corresponds to a point in $dS^3$. It follows that each geodesic plane of $\mathbb{H}^3$ corresponds to a point in $dS^3$.\\
Let $\Tilde{S}$ be a strictly convex surface of $\mathbb{H}^3$, and let $n: \Tilde{S} \to dS^3$ be the Gauss map, that is, for any $x \in \tilde{S}$, $n(x)$ is the unit normal vector on $T_{x}\tilde{S}$ pointing to the concave side. $n(\tilde{S}) = \tilde{S}^{*}$ is a strictly convex spacelike surface of $dS^3$, by spacelike we mean that the restriction of $\left<\cdot,\cdot\right>_{3,1}$ on $\Tilde{S}^{*}$ induces a Riemannian metric denoted by $I^{*}$, and by strictly convex we mean that its shape operator is positive definite. Moreover $(\Tilde{S}^{*},I^{*})$ is isometric to $(\Tilde{S},\III)$, we call $\Tilde{S}^{*}$ the dual surface of $\Tilde{S}$ (for more details see for example \cite[Section 1]{MR2208419}).\\
If $\Tilde{S}$ is the lift of a strictly convex surface $S$ embedded in a quasi-Fuchsian manifold $Q$ (see Section 2 for a definition of quasi-Fuchsian manifolds), then $\Tilde{S}^{*}$ determines a Cauchy surface $S^{*}$ inside a maximal global hyperbolic spatially compact de Sitter spacetime $Q^{*}$ (see \cite[Section 5]{MR1685590}). We call $S^{*}$ the dual surface of $S$.\\
Let $\tilde{S}$ be a convex surface in $\mathbb{H}^3$ which is not necessarily smoothly embedded, then we can define the dual tree of $\tilde{S}$ (for definitions and details we refer to \cite[Section 5]{MR1685590}). For any couple $(x,P_{x})$ where $x \in \tilde{S}$ and $P_{x}$ is a support plane on $\tilde{S}$ at $x$, we associate to $(x,P_x)$ the unit normal vector on $n_{P_{x}}$ pointing to the concave side. We denote the tree obtained in this way as $\tilde{S}^{*}$, It corresponds to the dual tree of $S$ ($\tilde{S^{*}}$ is the image of $(x,P_x)$ by $n_{.}$  when $x$ runs over $\Tilde{S}$ and $P_{x}$ runs over the support planes to $\Tilde{S}$ at $x$).\\

\section{Quasi-Fuchsian manifolds}
In this section, we define the objects we will study in this paper and give some key theorems concerning them.\\
We denote by $\mathbb{H}^3$ the hyperbolic space and we denote by $Isom^+(\mathbb{H}^3)$ the group of isometries of $\mathbb{H}^3$ that preserve orientation. A Kleinian group is a discrete sub-group of $Isom^+(\mathbb{H}^3)$. Recall that every  complete Riemannian manifold with a metric of constant sectional curvature equal to $-1$  is isometric to $\mathbb{H}^3/\Gamma$, where $\Gamma$ is a Kleinian group. For any Kleinian group $\Gamma$ and any point $p \in \mathbb{H}^3$, we define the orbit of $p$ under $\Gamma$ as  
$ \Gamma(p) := \left\{\gamma(p) \mid \gamma \in \Gamma \right\} $.  
The limit set of $\Gamma$, denoted by $\Lambda_{\Gamma}$, is the set of accumulation points of $\Gamma(p)$ in the boundary at infinity $\partial_{\infty} \mathbb{H}^3$. In other words,  
$ \Lambda_{\Gamma} = \overline{\Gamma(p)} \cap \partial_{\infty} \mathbb{H}^3 $,  
where the closure is taken in $ \mathbb{H}^3 \cup \partial_{\infty} \mathbb{H}^3 $. This definition does not depend on the choice of $p$. The complement of $\Lambda_ {\Gamma}$ in $\partial_{\infty}\mathbb{H}^3$ is referred to as the domain of discontinuity, which we denote by $\Omega_{\Gamma}$, that is:
$$\Omega_{\Gamma} := \partial_{\infty}\mathbb{H}^3 \setminus \Lambda_{\Gamma}.$$
The group $\Gamma$ acts properly discontinuously on $\Omega_{\Gamma}$ \cite[Proposition 8.2.3]{thurston2022geometry}.\\
Assume that $\Gamma$ is a finitely generated torsion free Kleinian group. We say that $\Gamma$ is a quasi-Fuchsian group if its limit set $\Lambda_{\Gamma}$ is a Jordan curve ( then a quasi-circle, see \cite[Proposition 8.7.2]{thurston2022geometry}) and each component of $\Omega_{\Gamma}$ is invariant under the action of $\Gamma$. In that case, $M_{\Gamma} = \mathbb{H}^3 / \Gamma$ is a hyperbolic manifold with conformal boundary $\Omega_{\Gamma} / \Gamma$ that consists of the disjoint union of two surfaces, we call such a hyperbolic manifold quasi-Fuchsian. If we add the assumption that $\Lambda_{\Gamma}$ is a geometric circle, we say that $\Gamma$ is Fuchsian. In the case when $\Gamma$ is a Fuchsian group we say that $M_{\Gamma}=\mathbb{H}^3/\Gamma$ is a Fuchsian manifold.\\
Assume that $\Gamma$ is quasi-Fuchsian, and let $\Omega^{\pm}_{\Gamma}$ be the components of $\Omega_{\Gamma}$. It  was shown by Marden in \cite{marden1974geometry} that $M_{\Gamma}=\mathbb{H}^3/\Gamma$ is diffeomorphic to $(\Omega^{+}_{\Gamma}/\Gamma) \times (0,1  )$ and that $\Bar{M}_{\Gamma}=(\mathbb{H}^3 \cup \Omega_{\Gamma})/\Gamma$ is diffeomorphic to $(\Omega^{+}_{\Gamma}/\Gamma) \times \left [0,1  \right ]$.  In this paper we assume that $\Gamma$ is a quasi-Fuchsian group that has no parabolic element, it yields that $\Omega^{+}_{\Gamma}/\Gamma$ is a closed surface of genus $g \geq 2$.\\

Let $M=S \times (0,1)$, and let $g$ be a Riemannian metric on $M$. We say that $(M,g)$ is quasi-Fuchsian if it is isometric to $\mathbb{H}^3 /\Gamma$ for some quasi-Fuchsian group $\Gamma$, and we say that $(M,g)$ is Fuchsian if $\Gamma$ is a Fuchsian group. We identify two metrics $g_1$ and $g_2$ if there is an isometry $\Phi : (M,g_1) \to (M,g_2)$ which is isotopic to the identity. We denote by $\mathcal{QF}(S)$ the equivalent classes of quasi-Fuchsian metrics on $M$ up to identification, and we denote by $\mathcal{F}(S)$ the subset of Fuchsian manifolds. Note that $\mathcal{QF}(S)$ has a topology induced from the set of representations $Hom(\pi_1(S),\mathbb{P}SL_{2}(\mathbb{C}))$. Through this paper when we mention a quasi-Fuchsian metric on $M$ we mean its isotopy class, that is, we identify any two metrics that are isometric to each other via an isometry isotopic to the identity.\\
Let $Q := \mathbb{H}^3 / \Gamma$ be a quasi-Fuchsian manifold. We define the ideal boundary $\partial_{\infty} Q$ of $Q$ as the disjoint union $\partial_{\infty}^{+} Q \sqcup \partial_{\infty}^{-} Q$, where $\partial_{\infty}^{+} Q$ (resp $\partial_{\infty}^{-} Q$) is identified with $\Omega^{+}_{\Gamma} / \Gamma$ (resp $\Omega^{-}_{\Gamma} / \Gamma$). Then, any quasi-Fuchsian metric $g$ on $M$ induces a conformal structure $c^+$ (resp $c^-$) on $S \times \{1\}$ (resp $S \times \{0\}$), such that $(S, c^+)$ (resp $(S, c^-)$) is identified with $\Omega^{+}_{\Gamma} / \Gamma$ (resp. $\Omega^{-}_{\Gamma} / \Gamma$).

Then we have a well defined map :
\begin{align*}
B : \mathcal{QF}(S) & \to \mathcal{T}(S) \times\mathcal{T}(S)\\
(M,g) & \mapsto  (c^+,c^-).
\end{align*}
A well-known theorem, known as Bers simultaneous uniformization theorem, is the following \cite{bams/1183523461}:  
\begin{theorem}\label{bers}
The map :  \begin{align*}
B : \mathcal{QF}(S) & \to \mathcal{T}(S) \times \mathcal{T}(S)\\
Q & \mapsto  (c^+,c^-)
\end{align*}
is a homeomorphism.
\end{theorem}
\subsection{Convex core} 
\begin{figure}
    \centering
    \includegraphics[scale=0.6]{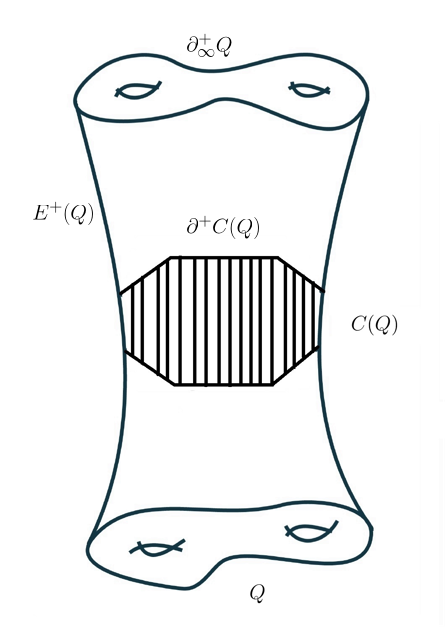}
\caption{The shadowed region is the convex core of this quasi-Fuchsian manifold.}
    \label{convex core-fig}
\end{figure}

For more details see \cite[Section 8.5]{thurston2022geometry}, \cite[Section II.1]{canary2006fundamentals}.
 Let $Q = (M,g) \in \mathcal{QF}(S)$, then Q is isometric to $\mathbb{H}^3/\Gamma$. Let $CH(\Lambda_{\Gamma})$ be the convex hull of $\Lambda_{\Gamma}$. We call the convex core of $Q$ the set $C(Q):= CH(\Lambda_{\Gamma})/\Gamma$ (we think of $C(Q)$ as a subset of $M$ after identification by isometry). The set $C(Q)$ is the smallest non-empty geodesically convex subset of $Q$ which contains every closed geodesic of $Q$. That is, if $C$ is a geodesically convex subset of $Q$ which contains every closed geodesic of $Q$, then $C(Q) \subset C$. Moreover the inclusion map $\iota : C(Q) \to \Bar{M}$ is a homotopy equivalence.\\
Except in the case when $Q$ is Fuchsian, the convex core $C(Q)$ has nonempty interior and its boundary $\partial C(Q)$ consists of the disjoint union of two surfaces $\partial^+ C(Q)$, and $\partial^-C(Q)$ and each of the two is homeomorphic to $S$. If $Q$ is Fuchsian, then $C(Q)$ is a totally geodesic surface homeomorphic to $S$.\\
Note that $M \setminus C(Q)$ consists of two connected components. We denote by $E^{+}(Q)$ the connected component whose closure in $\overline{M}$ contains $\partial^+_{\infty}Q$ and we denote by $E^{-}(Q)$ the connected component whose closure in $\overline{M}$ contains $\partial^-_{\infty}Q$. We assume by convention that $\partial^{+}C(Q) \subset \overline{E}^{+}(Q)$ and $\partial^{-}C(Q) \subset \overline{E}^{-}(Q)$. In particular, note that $\partial^{+}C(Q)$ (resp. $\partial^{-}C(Q)$) and $\partial^{+}_{\infty}Q$ (resp. $\partial^{-}_{\infty}Q$) belong to the same connected component of $\overline{M} \setminus C(Q)$ (see Figure \ref{convex core-fig}). It was shown by Thurston (see references above) that if $Q$ is not Fuchsian, then $\partial^+ C(Q)$ and $\partial^-C(Q)$ are pleated surfaces and $Q$ induces a path metric on each component which is hyperbolic, this gives two hyperbolic metrics $h^+$ and $h^-$.\\
Then we have a well defined map.
\begin{align*}
    T : \mathcal{QF}(S) & \to \mathcal{T}(S) \times \mathcal{T}(S)\\
    Q & \mapsto  (h^+,h^-).
\end{align*}
Thurston made the following conjecture,
\begin{conjecture}\label{thurston}(Thurston)
For any two hyperbolic metrics $h^{\pm}$ on $S$ (up to isotopy) there exists a unique quasi-Fuchsian manifold $Q$(up to isotopy) that induce $h^{\pm}$ on $\partial^{\pm} C(Q)$ (up to isotopy).
\end{conjecture}
The existence part of this statement is known since work of Labourie \cite{labourie1992metriques}, Epstein and Marden \cite{epstein1987analytical} and Sullivan \cite{sullivan2006travaux}. In other words, the conjecture says that the map $T$ is bijective, for the moment we only know that it is surjective.\\
The conformal structures at infinity and the induced metrics on the boundary of the convex core are related by the following theorem \cite[Section II.2]{canary2006fundamentals}.
\begin{theorem}(Sullivan, Epstein-Marden)\label{sullivan}
  There exists a universal constant $K$ such that $c^{+}$(respectively $c^{-}$) is $K$ quasi-conformal to $h^{+}$(respectively $h^{-}$).  
\end{theorem}
Combining Theorem \ref{sullivan} with Theorem \ref{bers}, we obtain the following corollary.
\begin{corollary}\label{comsull}
Let $(Q_{n})_{n \in \mathbb{N}} \subset \mathcal{QF}(S)$ be a sequence of quasi-Fuchsian manifolds such that $(T(Q_{n}))_{n \in \mathbb{N}}$ is contained in a compact subset of $\mathcal{T}(S) \times \mathcal{T}(S)$. Then $(Q_{n})_{n \in \mathbb{N}}$ has a convergent subsequence.  
\end{corollary}
\begin{proof}
  Recall that if $h$ and $h^{'}$ are two hyperbolic metrics on $S$ such that there is a $K$-quasi-conformal map which is homotopic to the identity between $(S,h)$ and $(S,h^{'})$, then by   \cite[Lemma 3.1]{wolpert1979length} we have the following inequality on the length spectrum 
  $$ \forall \alpha \in \mathcal{S}, \ \frac{\ell_{h^{'}}(\alpha)}{K} \leq \ell_{h}(\alpha) \leq K\ell_{h^{'}}(\alpha).$$
  Therefore, if $(T(Q_{n}))_{n \in \mathbb{N}}$ is contained in a compact subset of $\mathcal{T}(S) \times \mathcal{T}(S)$, Theorem \ref{sullivan} implies that $(B(Q_{n}))_{n \in \mathbb{N}}$ is also contained in a compact subset of $\mathcal{T}(S) \times \mathcal{T}(S)$.
  This implies that $(B(Q_{n}))_{n \in \mathbb{N}}$ has a convergent subsequence, and since $B$ is a homeomorphism by Theorem \ref{bers}, we deduce that $(Q_{n})_{n \in \mathbb{N}}$ has a convergent subsequence.
\end{proof}
As a consequence of Corollary \ref{comsull}, we get that each fiber of $T$ is compact.
\subsection{Bending laminations on the boundary of the convex core}
Let $Q = (M,g) \in \mathcal{QF}(S)$. Assume that $Q$ is not Fuchsian. Although the induced metrics on the boundary of the convex core $C(Q)$ are hyperbolic, the boundary surfaces of $C(Q)$ are not smoothly embedded, in particular they are not totally geodesic. If $Q$ is Fuchsian, $C(Q)$ turns to be a totally geodesic surface.\\
Thurston has noticed that $\partial^{\pm} C(Q)$ are convex pleated surfaces, and the pleating locus gives a measured lamination (see \cite{thurston2022geometry} section 8.5 and \cite{canary2006fundamentals} in section II.1). We have seen in Section 1.3, since $S$ is connected the pleating locus is a geodesic lamination, so $\partial^{+}C(Q)$ (resp $\partial^{-}C(Q)$) gives a geodesic lamination $\mu^{+}$(resp $\mu^{-}$) on $S$. Moreover, the amount of bending gives a transverse measure on each of $\mu^{\pm}$. The simplest case is when $\mu^{\pm}$ are rational, so we associate to each leaf an atomic weight which is equal to the exterior dihedral angle of the bending at that leaf. In the case when the lamination is not rational, it was shown in \cite[Section II.1.11]{canary2006fundamentals} that we can always approximate by rational ones. One remark is that any closed leaf of $\mu^{\pm}$ has an atomic weight strictly smaller than $\pi$.\\
It follows that we have a well defined map:
\begin{align*}
    L : \mathcal{QF}(S) & \to \mathcal{ML}(S) \times \mathcal{ML}(S)\\
    Q & \mapsto  (\mu^+,\mu^-).
\end{align*}
This map is clearly not surjective, but we know exactly what its image is.\\
Let $\mu^+$ and $\mu^-$ be two measured laminations on $S$, we say that they fill $S$ if there exists $\epsilon > 0$ such that for any simple closed curve $\alpha$ not homotopic to a point, $$i(\alpha,\mu^+)+i(\alpha,\mu^-) > \epsilon.$$
We define the set $\mathcal{ML}_{\pi}(S)$ to be the subset of $\mathcal{ML}(S)$ that consists of measured laminations in which every closed leaf has weight strictly smaller than $\pi$, and we define
$\mathcal{L}$ to be $$\mathcal{L} :=  \left\{(\mu^+,\mu^-) \in \mathcal{ML}_{\pi}(S) \times \mathcal{ML}_{\pi}(S)|  \mu^+ \ and \ \  \mu^- \textit{fill} \ S\right\}.$$

Thurston made the following conjectures.
\begin{conjecture}\label{bending}(Thurston)
Given any pair of measured laminations $\mu^+,\mu^- \in \mathcal{ML}_{\pi}(S)$ that fill $S$. There is a unique quasi-Fuchsian manifold $Q$ (up to isotopy), such that $\mu^{+}$ is the bending lamination of $\partial^{+}C(Q)$ and $\mu^{-}$ is the bending lamination of $\partial^{-}C(Q)$.
\end{conjecture}
Bonahon and Otal (see \cite{bonahon2004laminations}) proved the existence part, while Dular and Schlenker (see \cite{Dular-schlenker}) established the uniqueness.\\
In other words the conjecture says that the map 
\begin{align*}
 L: \mathcal{QF}(S) \setminus \mathcal{F}(S) \to & \mathcal{L}\\
 Q \mapsto & (\mu^{+},\mu^{-}).
\end{align*}
is bijective. Bonahon and Otal showed that the image of the map $L$ is exactly the set $\mathcal{L}$, while Dular and Schlenker proved that the map $L$ is injective.\\  
We also know that the maps $T$ and $L$ are continuous (see \cite{pjm/1102620682}). 

\subsection{Foliation by k-surfaces}
\begin{figure}
    \centering
    \includegraphics[scale=0.6]{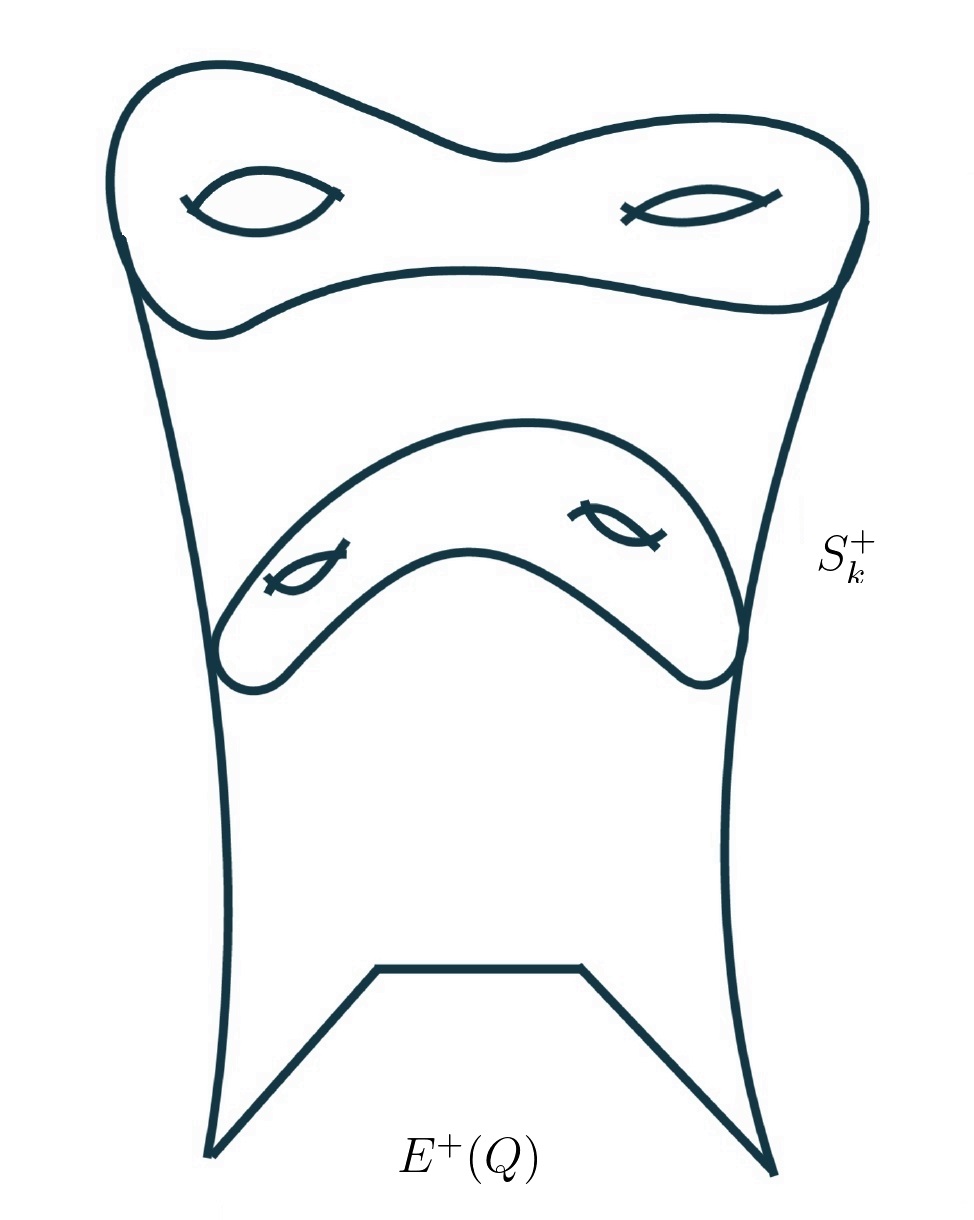}
\caption{The connected components of $Q \setminus C(Q)$ are foliated by k-surfaces.}
    \label{foliated-k-surfaces-fig}
\end{figure}

Let $Q \in \mathcal{QF}(S)$, and let $E^{+}(Q)$ and $E^{-}(Q)$ be the upper and  the lower connected components of $Q \setminus C(Q)$. It was shown by Labourie in \cite{labourie1991probleme} that $E^{+}(Q)$ (resp $E^{-}(Q)$) is foliated by surfaces $(S^+_{k})_{-1 < k < 0}$ (resp $(S^-_{k})_{-1 < k < 0}$) such that the induced metric on each $S^{\pm}_{k}$ has constant Gaussian curvature equal to $k$, and  each $S^{\pm}_{k}$ is homeomorphic to $S$ (see Figure \ref{foliated-k-surfaces-fig}). Each of $S^{\pm}_{k}$ is the unique surface in $E^{\pm}(Q)$ that has constant Gaussian curvature equal to $k$, and in which the projection $r^{\pm}_{k}:\partial^{\pm}_{\infty}Q \to S^{\pm}_{k}$ is a homotopy equivalence.\\
Moreover, for any $-1 < k < 0$ and $Q \in \mathcal{QF}(S)$, there is a unique convex subset $C_{k}(Q) \subset Q$, such that the inclusion map $\iota : C_{k}(Q) \to \Bar{Q}$ is a homotopy equivalence and $\partial C_{k}(Q) = S^+_{k} \cup S^-_{k}$.\\
From that, we deduce that any quasi-Fuchsian manifold induces two metrics $h^{+}_{k}$ and $h^{-}_{k}$ on $S$ of constant sectional curvature equal to $k$ by taking the induced metrics on $S^{+}_{k}$ and $S^{-}_{k}$ respectively. The quasi-Fuchsian manifold also induces two metrics $g^{+}_{k}$ and $g^{-}_{k}$ of curvature $\frac{k}{1+k}$ by taking the induced third fundamental forms on $S^{+}_{k}$ and $S^{-}_{k}$ respectively (recall that the Gaussian curvature of $(S,\III_{k})$ is equal to $\frac{k}{1+k}$, see Section 1.4).\\
This defines two new maps, 
\begin{align*}
 T_{k}: & \mathcal{QF}(S) \to \mathcal{T}(S) \times \mathcal{T}(S) \\  
 Q & \mapsto (-kh^+_{k},-kh^-_{k}),
\end{align*}
and 
\begin{align*}
 L_{k}: & \mathcal{QF}(S) \to  \mathcal{T}(S) \times \mathcal{T}(S) \\  
 Q & \mapsto  (\frac{-k}{1+k}g^+_{k},\frac{-k}{1+k}g^-_{k}).  
\end{align*}
Then there is a theorem which can be seen as a smooth version of Thurston's conjecture.
\begin{theorem}\label{smooth}
\noindent
\begin{itemize}
    \item Given any two metrics $h^{+}_{k}$ and $h^{-}_{k}$ on $S$ of curvature $k$, where $k \in ( -1,0 )$, there is a unique (up to isotopy) quasi-Fuchsian manifold that induces $h^{+}_{k}$ on $S^{+}_{k}$ and $h^{-}_{k}$ on $S^{-}_{k}$ as the first fundamental form (up to isotopy).
    \item Given any two metrics $g^{+}_{k}$ and $g^{-}_{k}$ on $S$ of curvature $\frac{k}{1+k}$, where $k \in (-1,0 )$, there is a unique (up to isotopy) quasi-Fuchsian manifold that induces $g^{+}_{k}$ on $S^{+}_{k}$ and $g^{-}_{k}$ on $S^{-}_{k}$ as the third fundamental form (up to isotopy).
\end{itemize}
\end{theorem}
The existence part of the first point of Theorem \ref{smooth} was proved by Labourie in \cite{labourie1992metriques}, and the uniqueness part was proved by Schlenker in \cite{MR2208419}. The second point of Theorem \ref{smooth} was proved by Schlenker in \cite{MR2208419}. In other words the maps $T_{k}$ and $L_{k}$ are bijective.\\
The maps $T_{k}$ and $L_{k}$ can be seen as an extension of the maps $T$ and $L$ when $k > -1$, this is because of the following theorem.
\begin{theorem}\label{extention}
 Let $(k_{n})_{n\in \mathbb{N}}$ be a decreasing sequence of negative numbers that converges to $-1$. Let $(Q_{n})_{n \in \mathbb{N}} \subset \mathcal{QF}(S)$ be a sequence of quasi-Fuchsian manifolds that converge to $Q_{\infty} \in \mathcal{QF}(S)$. We denote $T_{k_{n}}(Q_{n}) = (-k_{n}h^{+}_{k_{n}},-k_{n}h^{-}_{k_{n}})$ and $L_{k_{n}}(Q_{n}) = (-\frac{k_{n}}{1+k_{n}}g^+_{k_{n}},-\frac{k_{n}}{1+k_{n}}g^-_{k_{n}})$. Then the following assertions are true:
 \begin{enumerate}
     \item  $(h^+_{k_{n}},h^-_{k_{n}})$ converge to $T(Q_{\infty})$ in the length spectrum.
     \item $(g^+_{k_{n}},g^-_{k_{n}})$ converge to $L(Q_{\infty})$ in the length spectrum.
 \end{enumerate}
\end{theorem}
\begin{proof}
   The first point is a classical statement (see for example \cite[Proposition 6.6]{10.2140/gt.2013.17.157}).\\
The second point has been shown in details in \cite[Section 6]{10.2140/gt.2013.17.157}, for the reader's convenience we sketch the proof. Let's work on $S^{+}_{k_n}$ and denote it just by $S_{k_n}$. $S^{-}_{k_n}$ is treated in a similar way. As explained in Section 1.5, for each $n$, each $E^{+}(Q_{n})$ and $S_{k_n}$ determine a maximal globally hyperbolic spatially compact de Sitter spacetime  $Q^{*}_{_n}$ and a Cauchy surface $S^{*}_{k_n}$ which is dual to $S_{k_n}$. The foliation of $E^{+}(Q_n)$ by k-surfaces determines a foliation of $Q^{*}_{n}$ by surfaces $S^{*}_{k}$ of constant Gaussian curvature equal to $\frac{k}{1+k}$ on each leaf. Recall that $(S_{k_{n}},\III_{k_n})$ is isometric by the duality explained in Section 1.4 to $(S^{*}_{k_n},I^{*}_{k_n})$. Moreover, the initial singularity of $Q^{*}_{n}$ (which is a real tree), is dual to the bending lamination of $\partial^{+}C(Q_{\infty})$ (see \cite[Section 9]{scannell1996flat} for definitions and proofs).
As a corollary of \cite[Theorem 2.10]{AIF_2014__64_2_457_0}, the intrinsic metrics of the surfaces $S^{*}_{k}$ converge, with respect to the Gromov equivariant topology, to the real tree dual
of the bending lamination of $\partial^{+}C(Q_{\infty})$. It follows that the length spectrum of the surfaces $(S_{k_n}, \III_{k_n})$ (which are isometric to $(S^{*}_{k_n}, I^{*}_{k_n})$) converge to the length spectrum of the bending lamination on $\partial^{+}C(Q_{\infty})$.

\end{proof}

\subsection{The mixed boundary data of the convex core }
Let $Q \in \mathcal{QF}(S)$, and denote by $h^{-}_{k}$ and $g^{+}_{k}$ the metrics on $S$ coming from the first fundamental form of $S^{-}_{k}$ and the third fundamental form of $S^{+}_{k}$. This defines the map:
\begin{align*}
    M_{k}: & \mathcal{QF}(S) \to \mathcal{T}(S) \times \mathcal{T}(S)\\
    Q & \mapsto (-\frac{k}{1+k}g^{+}_{k},-kh^{-}_{k}),
\end{align*}
and also the map
\begin{align*}
    M_{\infty,k}: & \mathcal{QF}(S) \to \mathcal{T}(S) \times \mathcal{T}(S)\\
    Q & \mapsto (-\frac{k}{1+k}g^{+}_{k},c^{-}),
\end{align*}
where $c^{-}$ is the conformal structure on $\partial_{\infty}^-Q$.
The following theorem was shown by Chen and Schlenker \cite[Theorem 1.6]{chen2022geometric}.
\begin{theorem}\label{mix}
Let $h$ and $h^{*}$ respectively be a Riemannian metrics on $S$, and denote their curvatures by $k_{h}$ and $k_{h^*}$ respectively (which are not necessarily constant). We assume that $-1 < k_{h}$ and $k_{h^{*}} < 1$. Moreover, we assume that every contractible closed geodesic with respect to $h^{*}$ has length strictly bigger than $2\pi$. Let $c \in \mathcal{T}(S)$ be a conformal structure on $S$. Then, there exists a convex  hyperbolic metric $g$ on $M=S\times(0,1)$, the interior of $\overline{M}=S\times[0,1]$, such that:
\begin{itemize}
    \item $g$ induces the first fundamental form (respectively third fundamental form) on $S\times \left\{ 0\right\}$ which is isotopic to $h$ (respectively $h^{*}$).
    \item $g$ induces on $S \times \{1\}$ (which is the boundary at infinity of a quasi-Fuchsian manifold) a conformal structure isotopic to $c$.

\end{itemize}
Furthermore, the surface $S \times \{0\}$ is smoothly embedded and convex in a quasi-Fuchsian manifold $Q_M$, and $M$ is a convex submanifold of $Q_M$ with the same homotopy type.
\end{theorem}
In particular, the authors show that the maps $Mix_{k}$ and $Mix_{\infty,k}$ are bijective.\\
In this paper we will define two maps that are extensions of $Mix_{k}$ and $Mix_{\infty,k}$ to the boundary of the convex core. Let 
\begin{align*}
    Mix_{-1}: & \mathcal{QF}(S) \to \mathcal{ML}_{\pi}(S) \times \mathcal{T}(S)\\
    Q & \mapsto (\mu^+,h^-)
\end{align*}
Where $\mu^+$ is the bending lamination on $\partial^+ C(Q)$ and $h^-$ is the induced metric on $\partial^- C(Q)$.\\
Another map is :
\begin{align*}
Mix_{\infty,-1}: & \mathcal{QF}(S) \to \mathcal{ML}_{\pi}(S) \times \mathcal{T}(S)\\
    Q & \mapsto (\mu^+,c^-).
\end{align*}
Where $\mu^+$ is the bending lamination on $\partial^+ C(Q)$ and $c^-$ is the induced conformal structure on the lower boundary component at infinity $\partial^{-}_{\infty}Q$ of $Q$.\\
The main theorem of this paper is proving that the maps $Mix_{-1}$ and $Mix_{\infty,-1}$ are surjective. We prefer to state the theorem in the following way:
\begin{theorem}\label{Main}
 Let $ \mu \in \mathcal{ML}_{\pi}(S)$ and let $h \in \mathcal{T}(S)$, then:
 \begin{itemize}
     \item There exists a quasi-Fuchsian manifold $Q$ such that the bending lamination of $\partial^{+}C(Q)$ is $\mu$, and the induced metric on $\partial^{-}C(Q)$ is $h$.  
     \item There exists a quasi-Fuchsian manifold $Q$ such that the bending lamination of $\partial^{+}C(Q)$ is $\mu$, and the induced conformal structure on $\partial^{-}_{\infty}Q_{\infty}$ is $h$. 
 \end{itemize}
 
\end{theorem}
 Before proving the surjectivity of the maps $Mix_{\infty,-1}$ and $Mix_{-1}$, we will need to show that they are proper. So we give the following theorem : 
 \begin{theorem}\label{properness}
  The maps $Mix_{\infty,-1}$ and $Mix_{-1}$ are proper.
  \end{theorem}
  Note that the properness of these maps is equivalent.
  \begin{lemma}
    The map $Mix_{\infty,-1}$ is proper if and only if $Mix_{-1}$ is proper.   
  \end{lemma}
  \begin{proof}
   Let $(Q_{n})_{n \in \mathbb{N}} \subset \mathcal{QF}(S)$, then the sequence $Mix_{-1}(Q_{n}) = (L^{+}_{n},h^{-}_{n})$ lies on a compact subset if and only if $Mix_{\infty,-1}(Q_{n}) = (L^{+}_{n},c^{-}_{n})$ lies on a compact subset. Indeed, by Theorem \ref{sullivan} $h^{-}_{n}$ lies on a compact subset if and only if $c^{-}_{n}$ does.   
  \end{proof}

\section{Closing lemma}
In this section, we will show that the maps $Mix_{\infty,-1}$ and $Mix_{-1}$ are proper. To show that, we need to give a compactness statement that we will call the closing lemma. The name "Closing lemma" was originally given by Bonahon and Otal in their work on laminations \cite[Section 2]{bonahon2004laminations}, where the author proved the properness of the map that associates to a quasi-Fuchsian manifold the bending laminations of the boundary of the convex core. Before giving the main statement of the this section, we need to explain the notion of pleated annulus.
\subsection{Pleated annulus}
For more details see \cite[proof of Lemma 2.1]{bonahon1986bouts}. Let $M$ be a 3-dimensional hyperbolic manifold, let $\alpha$ be a closed geodesic and let $\alpha^*$ be a piece-wise geodesic closed curve which is homotopic to $\alpha$ by the homotopy $A: S^1 \times \left [0,1  \right ] \to M$. We want to make $Im(A)$ hyperbolically simplicial. Let $a_1,...,a_n$ be the points on $S^1 \times \left\{ 1\right\}$, in this order, that are sent by $A$ to the  vertices of $\alpha^{*}$. Choose $b_1,...,b_n$ on $S^1 \times \left\{ 0\right\}$, in this order. We triangulate the annulus $S^1 \times \left [ 0,1 \right ]$ by joining $a_{i}$ to $b_{i}$ and $b_{i+1}$ ($b_{n+1} = b_{1}$). We can assume that each triangle is sent by $A$ to a totally geodesic triangle in $M$. Then the hyperbolic metric on $M$ induces a path metric on $Im(A)$ which is hyperbolic, and makes $Im(A)$ into a hyperbolic surface with piece-wise geodesic boundary.\\
An important fact for us is that by the Gauss-Bonnet Theorem the area of $Im(A)$ is equal to the sum of the exterior angles of $\alpha^*$.\\
In particular, if $Q \in \mathcal{QF}(S)$ and $\alpha^*$ lies on $\partial C(Q)$, then the area of the annulus $Im(A)$ is smaller than $i(\alpha,B)$, where $B$ is the bending lamination of the component of $\partial C(Q)$ on which $\alpha^*$ lies. Indeed, this is true because the sum of the exterior angles of $\alpha^*$ is smaller than the sum of the exterior dihedral angles of $\partial C(Q)$. This holds because each exterior angle of $\alpha^*$ corresponds to a corner where two consecutive geodesic segments meet. If $\alpha^*$ lies on the pleated surface $\partial C(Q)$, then at each such corner, the turn of $\alpha^*$ is controlled by the dihedral angle between the two faces of $\partial C(Q)$ that meet along the bending line crossed by $\alpha^*$. Since $\alpha^*$ follows the surface but does not necessarily bend as much as the surface itself, the sum of its exterior angles is strictly less than the total sum of the dihedral bending angles along the bending lamination crossed by $\alpha^*$ (see \cite[Lemma 3.2]{MR1193599}).\\
If $\alpha^*$ is not finitely bent, then we approximate it by finitely piecewise geodesic closed curves (see \cite[Chapter II]{canary2006fundamentals}).

\subsection{Closing lemma}
the following Lemma, that we call the closing lemma, is the main statement of this section.

\begin{lemma}\label{closing}
 Let $(Q_{n})_{n \in \mathbb{N}} \subset \mathcal{QF}(S) $ be a sequence of quasi-Fuchsian manifolds, denote by $h^{-}_n$ the induced metric on $\partial^{-}C(Q_n)$, and by $B^{+}_n$ the bending lamination of $\partial^{+}C(Q_n)$. Assume that: 
 \begin{itemize}
     \item $(h^{-}_n)_{n \in \mathbb{N}}$ belongs to a compact subset of $\mathcal{T}(S)$. 
     \item $(B^{+}_n)_{n \in \mathbb{N}}$ converges in the length spectrum to a measured lamination $B_{\infty}$, and every closed leaf of $B_{\infty}$ has weight strictly smaller than $\pi$.
 \end{itemize}
 Then up to extracting a subsequence, the sequence $(Q_{n})_{n \in \mathbb{N}}$ converges to a quasi-Fuchsian manifold $Q_{\infty}$.
\end{lemma}

Note that the statement of Lemma \ref{closing} means exactly that the map $Mix_{-1}$ (and therefore $Mix_{\infty,-1}$ ) is proper. 
The proof of this lemma is mainly based on the arguments given by \cite{lecuire2014convex} in Section 3. Nevertheless there are some differences between the two proofs, that's why we prefer to give the full arguments here.\\
Recall that by Corollary \ref{comsull} $(Q_n)_{n \in \mathbb{N}}$ converges (up to extracting a subsequence) if and only if both of $h^{+}_{n}$ and $h^{-}_{n}$, the induced metrics on $\partial^{+}C(Q_n)$ and $\partial^{-}C(Q_n)$ respectively, belong to a compact subset of $\mathcal{T}(S)$. Then the proof of Lemma \ref{closing} is based on showing for any $\alpha \in \mathcal{S}$, the sequence of lengths $(\ell_{h^{+}_{n}}(\alpha))_{n \in \mathbb{N}}$ is bounded  (recall that $\mathcal{S}$ is the set of free homotopy classes of simple closed curves not homotopic  to a point on $S$).\\
In what follows we will show that if $\ell_{h^{+}_{n}}(\alpha) \to \infty$, then $B_{\infty}$ must have a leaf of weight bigger than or equal to $\pi$. The fact that $(h^{-}_{n})_{n \in \mathbb{N}}$ belongs to a compact subset of $\mathcal{T}(S)$ plays a crucial role.\\

The previous discussion leads us to show the following lemma. 
\begin{lemma}\label{subclos}Let $(Q_{n})_{n \in \mathbb{N}} \subset \mathcal{QF}(S) $ be a sequence of quasi-Fuchsian manifolds, denote by $h^{-}_n$ (resp $h^{+}_n$) the induced metric on $\partial^{-}C(Q_n)$ (resp $\partial^{+}C(Q_n)$), and denote by $B^{+}_n$  the bending lamination of $\partial^{+}C(Q_n)$. If we assume that :
\begin{itemize}
    \item$\exists$ $\alpha \in \mathcal{S}$ such that $\ell_{h^{+}_{n}}(\alpha) \to \infty$.
    \item $(h^{-}_{n})_{n \in \mathbb{N}}$ is contained in a compact subset of $\mathcal{T}(S)$.
    \item $(B^{+}_{n})_{n \in \mathbb{N}}$ converge in the length spectrum to a measured lamination $B_{\infty}$.
\end{itemize}
Then $B_{\infty}$ has at least one closed leaf of weight bigger than or equal to $\pi$.
\end{lemma}
The key technical lemma for showing Lemma \ref{subclos} is the following : 
\begin{lemma}\label{pi}
  Under the hypothesis of Lemma \ref{subclos}, there is a sequence of geodesic arcs $(\kappa_{n})_{n \in \mathbb{N}}$ (each $\kappa_{n}$ is a geodesic arc of $Q_{n}$), such that each $\kappa_{n}$ has endpoints in $\partial^{+}C(Q_n)$, $\ell_{m_n}(\kappa_{n}) \to 0$, and $\ell_{h^{+}_{n}}(\kappa'_{n}) \to \infty$, where $\kappa'_{n}$ is the geodesic arc in $\partial^{+}C(Q_n)$ which is homotopic to $\kappa_{n}$ relative to its endpoints.\\
  Moreover, if we denote by $\alpha_{n}$ the unique simple closed geodesic representative of $\alpha$ in $\partial^{+}C(Q_n)$, then we can choose the endpoints of $\kappa_{n}$ on $\alpha_{n}$, and $\kappa'_{n}$ to be a subarc of $\alpha_{n}$.  
\end{lemma}
\begin{proof}
  Let $\alpha \in \mathcal{S}$ be such that $\ell_{h^{+}_n}(\alpha) \to \infty$, we denote by $\alpha_{n}$ the unique simple closed geodesic in $\partial^{+}C(Q_n)$ that belongs to $\alpha$, and we denote by $\alpha^{*}_{n}$ the unique simple closed geodesic of $Q_n$ which is freely homotopic to $\alpha_{n}$. Since $\ell_{h^{-}_n}(\alpha) \geq \ell_{m_n}(\alpha^{*}_{n})$ (because $\alpha^{*}_{n}$ minimizes the lengths on its free homotopy class. Here, by $m_n$ we mean the quasi-Fuchsian metric on $Q_n$), and since the sequence $(\ell_{h^{-}_n}(\alpha))_{n \in \mathbb{N}}$ is bounded (by hypothesis), we deduce that the sequence $(\ell_{m_n}(\alpha^{*}_{n})_{n \in \mathbb{N}})$ is also bounded by some $L > 0$.\\
  To continue the proof, we need to show the following technical sublemma : 
  \begin{sublemma}\label{subarc}
Up to extracting a subsequence, there is a sequence $(\eta_{n})_{n \in \mathbb{N}}$ such that each $\eta_{n}$ is a subarc of $\alpha_{n}$ ($\alpha_{n}$ is the geodesic in $\partial^{+}C(Q_n)$ with homotopy class $\alpha$), $\ell_{m_{n}}(\eta_{n}) \to \infty$, and $i(\eta_{n},B^{+}_{n}) \to 0$.    
  \end{sublemma}
  \begin{proof}
We argue by contradiction. If the lemma were false, then there exists $r > 0$ such that for every sequence of subarcs $\eta_n \subset \alpha_n$ with $\ell_{m_n}(\eta_n) \to \infty$, we have  $\displaystyle\liminf\limits_{n \to \infty} i(\eta_n, B_n^+) > r$.\\  
Then, for any $N \in \mathbb{N}$, we can partition $\alpha_n$ into $N$ subarcs, each of which has length tending to infinity. This implies that $ i(\alpha_n, B_n^+) > Nr$ for any $N \in \mathbb{N}$, so  
$\displaystyle\lim\limits_{n \to \infty} i(\alpha_n, B_n^+) = \infty$,  
which contradicts the hypothesis that $i(\alpha_n, B_n^+) \to i(\alpha, B_\infty)$. This completes the proof.

  \end{proof}
  Let $\eta'_{n}$ be the geodesic arc in $Q_n$ which is homotopic to $\eta_{n}$ relative to its endpoints. By \cite[Lemma A.1]{MR2207784}, we know that the Hausdorff distance between $\eta_{n}$ and $\eta'_{n}$ converges to $0$ and that $\ell_{m_{n}}(\eta'_{n}) \to \infty$.\\
  Let $A_n$ be the annulus bounded by $\alpha_{n}$ and $\alpha^{*}_{n}$ as explained in Section 3.1. The area of $A_n$ is bounded by $i(\alpha_{n},B^{+}_{n})$. Let $A'_{n}$ be the same annulus as $A_n$, except that we replace $\eta_{n}$ by $\eta'_{n}$, that is, $A'_{n}$ is a pleated annulus bounded by $\alpha^{*}_{n}$ and $\alpha'_n := (\alpha_{n}\setminus \eta_{n}) \cup \eta'_{n}$ (see Figure \ref{schema-of-proof-fig}). Note that the sums of the angles of  $\alpha'_n$ are bounded independently on $n$ (because the sums of the angles of  $\alpha_n$ are bounded), this ensures that the areas of $A'_{n}$ are bounded independently on $n$ by a constant $K >0$.\\

Let $a_n,b_n \in \eta'_{n}$ be the endpoints of $\eta'_{n}$, and let $a'_n,b'_n \in 
\eta'_{n}$ such that $d_{\eta'_{n}}(a_n,a'_{n}) = \frac{\ell_{m_{n}}(\eta'_{n})}{3}$ and $d_{\eta'_{n}}(b_n,b'_{n}) = \frac{\ell_{m_{n}}(\eta'_{n})}{3}$, where $d_{\eta'_{n}}$ is the distance induced from the metric $m_{n}$ of $Q_n$ restricted to the arc $\eta'_{n}$. Denote by $\nu_{n}$ the subarc of $\eta'_{n}$ that has $a'_{n}$ and $b'_{n}$ as endpoints, then $\ell_{m_{n}}(\nu_{n})= \frac{\ell_{m_{n}}(\eta'_{n})}{3}$ (that is $\nu_{n}$ is the second third segment of $\eta'_{n}$).\\

Let $\delta_{n} := \sinh^{-1}(\frac{9K}{\ell_{m_{n}}(\eta'_{n})} + \frac{1}{n})$. Denote by $E_{n}$ the set of segments in $A'_{n}$ that are orthogonal on $\nu_{n}$, have one endpoint in $\nu_{n}$, and either have length equal to $\delta_{n}$ and exactly one endpoint on $\partial A'_{n}$, or have length strictly smaller than $\delta_{n}$ and their second endpoints belong to $\partial A'_{n}$ (the first endpoint belongs to $\partial A'_{n}$ by definition). We denote by $D_{n}$ the subset of $E_{n}$ that consists of segments having length equal to $\delta_{n}$, and we denote by $Z_{n}$ the set of their endpoints on $\nu_{n}$.\\

\begin{sublemma}\label{help}
  If  for any $n \in \mathbb{N}$ there exists a segment $\xi_{n} \in E_{n} \setminus D_{n}$ that has both endpoints on $\alpha'_n$, then Lemma \ref{pi} holds. 
\end{sublemma}
\begin{proof}
   Assume the existence of a segment $\xi_{n} \in E_{n} \setminus D_{n}$ that has both endpoints on $\alpha'_n$. Because $\xi_n$ is orthogonal to $\eta'_{n}$, it has one endpoint on $\nu_{n}$ and the other endpoint will be on $\alpha'_{n}\setminus \eta'_{n} = \alpha_{n} \setminus \eta_{n}$. Denote these two endpoints by $x_n$ and $y_{n}$ respectively ($x_n \in \nu_{n}$ and $y_n \in \alpha'_{n} \setminus \eta'_{n}$). Recall that the Hausdorff distance between $\eta_{n}$ and $\eta'_{n}$ goes to $0$. Because $x_n \in \nu_{n} \subset \eta'_{n}$, we can find a sequence of points $z_{n} \in \eta_{n}$ such that $\lim\limits_{n \to \infty} d_{m_{n}}(z_n,x_n) = 0$. Then $\lim\limits_{n \to \infty}d_{m_n}(z_n,y_n) = 0$, it yields that there exists a geodesic arc $\kappa_{n}$ in $Q_n$ joining $z_{n}$ and $y_{n}$ such that $\lim\limits_{n \to \infty}\ell_{m_{n}}(\kappa_{n})=0$. On the other hand, we have that $\kappa'_{n}$, the subarc of $\alpha_{n}$ that joins $z_{n}$ and $y_{n}$ which is homotopic to $\kappa_{n}$, relative to its end points, satisfies $\ell_{m_n}(\kappa'_{n}) \geq \frac{\ell_{m_{n}}(\eta'_{n})}{3} - d_{m_n}(x_n,z_n)$. Also note that $\kappa'_{n}$ is a geodesic arc of $\partial^+ C(Q_n)$. Moreover, note that $\kappa_{n}$ and $\kappa'_{n}$ can be chosen 
   in a way that the endpoints of $\kappa_{n}$ are on $\alpha_{n}$ and $\kappa'_{n}$ is a subarc of $\alpha_{n}$.\\ We conclude that $\kappa_{n}$ and $\kappa'_{n}$ satisfy the statement of Lemma \ref{pi}.\\
   
\end{proof}

\begin{figure}
    \centering
    \includegraphics[scale=0.6]{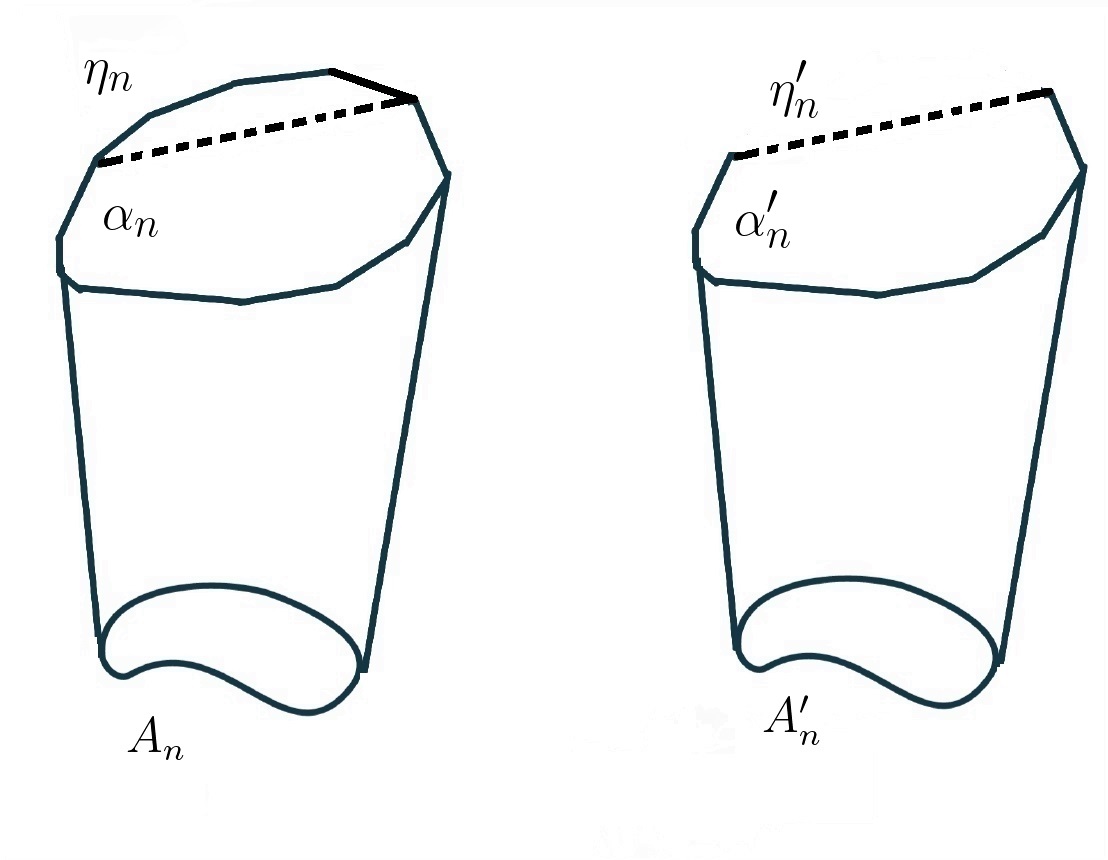}
\caption{We get the annulus $A'_{n}$ from the annulus $A_{n}$ by replacing the piece-wise geodesic arc $\eta_{n}$ by the geodesic arc $\eta'_{n}$.}
    \label{schema-of-proof-fig}
\end{figure}

In what follows we will use the fact that the lengths $\ell_{m_{n}}(\alpha^{*}_{n})$ are bounded to show that we are always in the case of Sublemma \ref{help}. We argue by contradiction, so we assume that every segment $\xi_n \in E_n$ that has length strictly smaller than $\delta_n$ does not intersect $\alpha'_n$. Since, by definition of $E_n$, the segment must intersect the boundary (because if the segment has length smaller than $\delta_n$, it must intersect $\partial A'_n$), it must intersect $\alpha^*_n$. Let $\mathcal{D}_{n} = \cup_{d \in D_n}d$, note that $\mathcal{D}_{n}$ is a hyperbolic strip (maybe disconnected), so $\area(\mathcal{D}_{n}) = \ell_{\eta'_{n}}(Z_{n}).\sinh(\delta_{n})$ (where $\ell_{\eta'_{n}}$ is the length measure induced on $\eta'_{n}$ by $m_{n}$), and since it is a subset of $A'_{n}$ and the $(\area(A'_{n}))_{n \in \mathbb{N}}$ are bounded by some constant $K$, it follows that,
$$K \geq \area(A'_{n}) \geq \area(\mathcal{D}_{n}) = \ell_{\eta'_{n}}(Z_{n}) \cdot \sinh(\delta_{n}) > \ell_{\eta'_{n}}(Z_{n}) \cdot \frac{9K}{\ell_{m_{n}}(\eta'_{n})}.$$
Then $\frac{\ell_{m_{n}}(\eta'_{n})}{9} > \ell_{\eta'_{n}}(Z_{n})$. In particular we can find two segments $\xi_{n},\xi'_{n} \in E_{n}\setminus D_{n}$ with endpoints $\xi_{n}(0),\xi'_{n}(0) \in \nu_{n}$, such that $\nu'_{n}$, the subarc of $\nu_{n}$ with endpoints  $\xi_{n}(0)$ and $\xi'_{n}(0)$, satisfies $\ell_{\eta'_{n}}(\nu'_{n}) \geq \frac{\ell_{m_{n}}(\eta'_{n})}{9}$ (for example take $\xi_{n}$ and $\xi'_{n}$ such that $\xi_{n}(0)$ is in the first third of $\nu_{n}$ and $\xi'_{n}(0)$ is in the last third of $\nu_{n}$). Let $\Bar{\alpha}_{n}$ be the subarc of $\alpha^{*}_{n}$ with endpoints $\xi_{n}(1)$ and $\xi'_{n}(1)$ and such that the arc $\zeta_{n}=\xi_{n} \cup \Bar{\alpha}_{n} \cup \xi'_{n}$ is homotopic to $\nu'_{n}$ relatively end points. This yields a contradiction, because on one hand the lengths $\ell_{m_{n}}(\kappa_{n})$ are bounded, and on the other hand $\ell_{m_{n}}(\nu'_{n}) \to \infty$, this can not happen because $\nu'_{n}$ is a geodesic arc so we must have $\ell_{m_{n}}(\nu'_{n}) \leq \ell_{m_{n}}(\kappa_{n}) $.\\
It follows that the hypothesis of Sublemma \ref{help} holds, and then Lemma \ref{pi} is true.
\end{proof}
The rest of the proof of Lemma \ref{subclos} follows from Claim 3.6, Claim 3.8 and Claim 3.9 from \cite{lecuire2014convex}. Nevertheless, we provide a proof for the reader's convenience. The proof that we provide is a bit simpler because we are dealing with quasi-Fuchsian manifolds without particles. \\
It will happen that we will use $i(\kappa,\lambda) = \int_{\kappa}d\lambda_{\kappa}$ as defined in Subsection 1.2 instead of taking the map $\mathcal{I}$ (see Section 1.2). When we use it, we will point it out.
\begin{lemma}\label{3.6}
 Let $\kappa_{n}$ and $\kappa'_{n}$ be the arcs constructed in Lemma \ref{pi}, and let $y_n$ and $z_n$ be the endpoints of $\kappa_{n}$ (therefore the endpoints of $\kappa'_{n}$). Let $l_n$ be a simple closed curve on $\partial^{+}C(Q_n)$ which is based at $y_n$ and which is geodesic in $\partial^{+}C(Q_n)$ (except at $y_n$). Let $f_{n}$ be the simple closed curve based at $z_n$ which is geodesic in $\partial^{+}C(Q_n)$ (except at $z_n$) and freely homotopic to $l_n$. Choose $(l_n)_{n \in \mathbb{N}}$ such that they have a bounded length (independently on $n$). Then $\lim\limits_{n \to \infty}i(l_{n},B^{+}_{n}) =0$ and $\lim\limits_{n \to \infty}i(f_{n},B^{+}_{n}) = 0$, where here we mean by $i(l_n,B^{+}_{n})$ and $i(f_{n},B^{+}_{n})$ the weight associated to $l_n$ and $f_n$ as transverse arcs (instead of taking the weight associated to the simple closed geodesic that represents them).    
\end{lemma}
\begin{proof}
 Assume that $Q_n \approx \mathbb{H}^3/G_n$, where $G_n$ is a quasi-Fuchsian group. Let $\Tilde{y}_n$ and $\Tilde{z}_n$ be lifts of $y_n$ and $z_n$ respectively, let $\Tilde{\kappa}_n$ be a lift of $\kappa_{n}$ and $\Tilde{\kappa}'_{n}$ be a lift of $\kappa'_n$, such that both of $\Tilde{\kappa}_n$ and $\Tilde{\kappa}'_{n}$ have $\Tilde{y}_n$ and $\Tilde{z}_n$ as endpoints, and let $\partial^{+}\Tilde{C}(Q_n)$ be the lift of $\partial^{+}C(Q_n)$. Up to composition by isometries of $\mathbb{H}^3$, we can assume that $\Tilde{y}_n$ is constant independently on $n$. Up to moving $\Tilde{y}_n$ and $\Tilde{z}_n$ slightly we can assume that they are away from the pleating locus of $\partial^{+}\Tilde{C}(Q_n)$.\\
 Let $P(\Tilde{y}_n)$ be the support plane at $\Tilde{y}_n$ and let $P(\Tilde{z}_n)$ be the support plane at $\Tilde{z}_n$. The fact that $\ell_{\partial^{+}C(Q_n)}(\kappa'_{n}) \to \infty$ implies that either $P(\Tilde{y}_n) \cap P(\Tilde{z}_n)$ is empty or it diverges to the boundary at infinity (that is, for any $K$, a compact subset of $\mathbb{H}^3$, the set $\left\{n \in \mathbb{N}| \  P(\Tilde{y}_n) \cap P(\Tilde{z}_n) \cap K \neq \emptyset \ \right\}$ is finite). And since $\ell_{m_n}(\Tilde{\kappa}_n) \to 0$, we deduce that $P(\Tilde{y}_n)$ and $P(\Tilde{z}_n)$ converge (up to extracting a subsequence) to the same plane, call it $P_{\infty}$ (see Figure \ref{figure-lemma 3.6}). We denote by $E(\Tilde{y}_n)$ (respectively $E(\Tilde{z}_n)$)  the half space determined by $P(\Tilde{y}_n)$ (respectively $P(\Tilde{z}_n)$) and contains $\partial^{+}\Tilde{C}(Q_n)$. Let $\Tilde{l}_n$ be a lift of  $l_n$ with endpoints $\Tilde{y}_n$ and $\Tilde{y}'_{n}$, let $g_n \in G_n$ such that $\Tilde{y}'_{n} = g_n \Tilde{y}_n$. Let $\Tilde{f}_n$ be a lift of  $f_n$ with endpoints $\Tilde{z}_n$ and $\Tilde{z}'_{n}$ such that $\Tilde{z}'_{n} = g_n\Tilde{z}_n$.\\
 \begin{figure}
     \centering
     \includegraphics[width=0.5\linewidth]{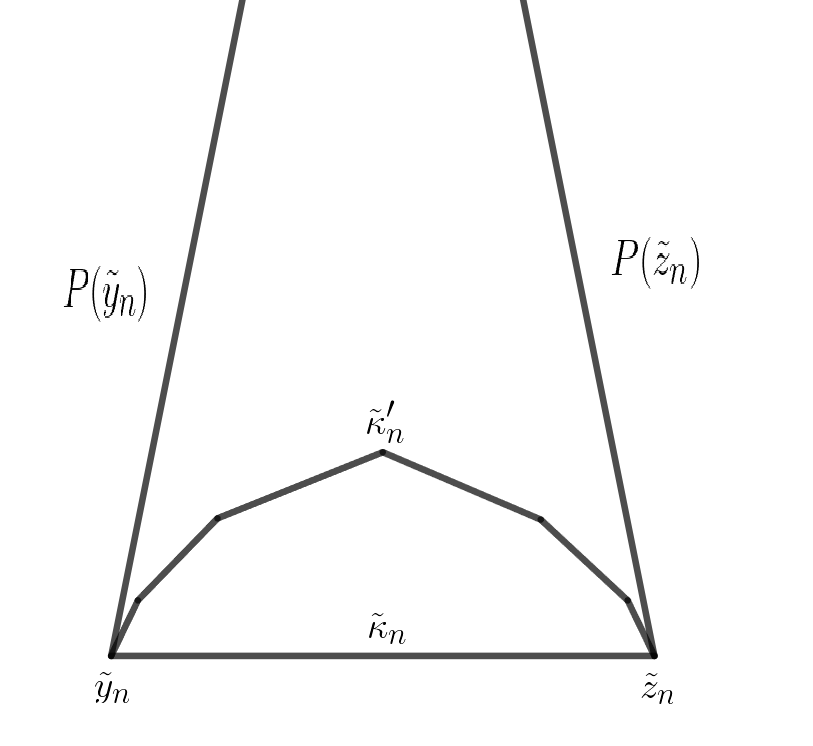}
     \caption{The fact that $\ell_{m_n}(\tilde{\kappa}_n) \to 0$ implies that the planes $P(\tilde{y}_n)$ and $P(\tilde{z}_n)$ converge (up to extracting a subsequence) to the same plane. Moreover, since $\ell_{m_n}(\tilde{\kappa}'_n) \to \infty$, we deduce that $P(\tilde{y}_n) \cap P(\tilde{z}_n)$ is either empty or diverges to the boundary at infinity.}
     \label{figure-lemma 3.6}
 \end{figure}
Let $P(\Tilde{y}'_{n})$ be the support plane at $\Tilde{y}'_{n}$ ( $\Tilde{y}'_{n}$ is away from the pleating locus, because the pleating locus is preserved by the action of the group $G_n$). Note that $d_{\partial^{+}\Tilde{C}(Q_n)}(\Tilde{y}'_{n},\Tilde{z}_n) \to \infty$ (because the lengths of $(\ell_{\partial^{+} C(Q_n)}(l_n))_{n \in \mathbb{N}}$ are bounded independently on $n$). Then $P(\Tilde{y}'_n) \cap P(\Tilde{z}_n)$ is empty or diverge to the boundary at infinity. Since $\Tilde{y}_n$ is fixed and since $\Tilde{y}'_{n}$ is at a bounded distance from $\Tilde{y}_n$ it follows that up to extracting a subsequence $P(\Tilde{y}'_n)$ converge to $P_{\infty}$ (Indeed, since $E(\tilde{y}_n) \cap E(\Tilde{z}_n)$ converge to  an empty set, it follows that if $P(\Tilde{y}'_n)$ does not converge to the same plane as $P(\Tilde{z}_n)$, then $P(\Tilde{y}'_n) \cap P(\Tilde{z}_n)$ will not be empty and will not go to the boundary at infinity). Now it follows that both of $P(\Tilde{y}_n)$ and $P(\Tilde{y}'_n)$ converge to $P_{\infty}$, and since $i(l_n,B^{+}_{n})$ is dominated by the dihedral angle of the intersection between $P(\Tilde{y}_n)$ and $P(\Tilde{y}'_n)$ \cite[Section II.1.10]{canary2006fundamentals}, we deduce that $i(l_n,B^{+}_{n}) \to 0$.\\
Finally, recall that $\Tilde{y}'_{n} = g_n\Tilde{y}_{n}$ and $\Tilde{z}'_{n} = g_n\Tilde{z}_{n}$, then $d_{\mathbb{H}^3}(\Tilde{y}'_{n},\Tilde{z}'_{n}) = d_{\mathbb{H}^3}(\Tilde{y}_{n},\Tilde{z}_{n}) $ and $d_{\partial^{+}\Tilde{C}(Q_n)}(\Tilde{y}'_{n},\Tilde{z}'_{n}) = d_{\partial^{+}\Tilde{C}(Q_n)}(\Tilde{y}_{n},\Tilde{z}_{n}) $. So for the same reasons as above $P(\Tilde{z}'_{n})$ converge up to extracting a subsequence to the same limit as $P(\Tilde{y}'_{n})$, therefore to the same limit as $P(\Tilde{z}_{n})$, therefore we deduce that $i(f_n,B^{+}_{n}) \to 0$.

\end{proof}

\begin{lemma}
  Let $\kappa_{n}$ and $\kappa'_{n}$ be the arcs constructed in Lemma \ref{pi}, then $\liminf i(\kappa'_{n},B^{+}_{n}) \geq \pi$.  
\end{lemma}
\begin{proof}
Up to approximating $\kappa'_n$ by piecewise geodesic segments, we may assume that the curve $\kappa_n \cup \kappa'_n$ is a skew polygon. Note that this curve bounds a disk in $C(Q_n)$. Let $\tilde{\kappa}_n$ and $\tilde{\kappa}'_n$ be lifts of $\kappa_n$ and $\kappa'_n$ respectively so that $\tilde{\kappa}_n \cup \tilde{\kappa}'_n$ is a lift of $\kappa_n \cup \kappa'_n$. Denote by $\tilde{y}_n$ and $\tilde{z}_n$ the endpoints of $\tilde{\kappa}_n$, and let $C_n$ be the geodesic cone with vertex $\tilde{y}_n$ spanning $\tilde{\kappa}_n \cup \tilde{\kappa}'_n$. In other words, $C_n$ is the union of hyperbolic triangles formed by connecting $\tilde{y}_n$ to each edge of the skew polygon $\tilde{\kappa}_n \cup \tilde{\kappa}'_n$ by geodesic segments.

Since $\tilde{\kappa}_n \cup \tilde{\kappa}'_n$ is piecewise geodesic, it follows that $C_n$ is a finite union of hyperbolic triangles and, hence, the induced metric on $C_n$ is hyperbolic. Moreover, since $\ell_{\mathbb{H}^3}(\tilde{\kappa}_n) \to 0$, the support planes at $\tilde{z}_n$ and $\tilde{y}_n$ intersect and converge to the same plane. Thus, the sum of the interior angles at $\tilde{y}_n$ and $\tilde{z}_n$ approaches $\pi$.

Recall that by the Gauss-Bonnet theorem we have :
\[
\sum_i \theta_i^n = \mathrm{Area}(C_n) + 2\pi,
\]
where the $\theta_i^n$ are the exterior angles at the vertices of $C_n$, and we denote by $\theta_{\tilde{z}}^n$ and $\theta_{\tilde{y}}^n$ the exterior angles at $\tilde{z}_n$ and $\tilde{y}_n$, respectively. From the discussion above, we have :
\[
\theta_{\tilde{z}}^n + \theta_{\tilde{y}}^n \to \pi.
\]
Therefore, applying the Gauss--Bonnet formula yields : 
\[
\liminf_{n\to\infty} i(\kappa'_n, B^+_n) \ge 2\pi - \limsup_{n\to\infty} (\theta_{\tilde{z}}^n + \theta_{\tilde{y}}^n) \ge \pi.
\]
\end{proof}

Now we are ready to give a proof of Lemma \ref{subclos}.
\begin{proof}
 We identify $Q_n \cong \mathbb{H}^3/G_n$. Let $\kappa_{n}$,$\kappa'_{n}$,$l_n$ and $f_n$ as defined in Lemma \ref{3.6}. Let $\Tilde{\kappa}_n$ be a lift of $\kappa_{n}$ with endpoints $\Tilde{y}_n$ and $\Tilde{z}_n$, let $\Tilde{\kappa}'_{n}$ be the lift of $\kappa'_n$ with endpoints $\Tilde{y}_n$ and $\Tilde{z}_n$. Take $\Tilde{l}_n$  (respectively $\Tilde{f}_n$) to be a connected component of the preimage of $l_n$ (respectively of $f_n$) such that $\Tilde{l}_n$ and $\Tilde{f}_n$ are connected by $\Tilde{\kappa}'_n$. The curves $\Tilde{l}_n$ and $\Tilde{f}_n$ are disjoint broken geodesics that bound an infinite band $B_n$ in $\partial^{+}\Tilde{C}(Q_n)$ (where $\partial^{+}\Tilde{C}(Q_n)$ is the lift of $\partial^{+}C(Q_n)$). Let $\left<g_n \right> \subset G_n$ be the subgroup generated by $g_n$, such that for all $n$, the action of $g_n$ leaves the elements $\Tilde{l}_{n}$ invariant. Since $f_n$ is homotopic to $l_n$ it follows that $\Tilde{f}_n$ is also invariant under the action of $\left<g_n \right>$. Let $c_n$ be the closed simple geodesic in the homotopy class of $l_n$ and $f_n$ and let $\Tilde{c}_n$ be the component of the preimage of $c_n$ that has same endpoints as $\Tilde{l}_n$ and $\Tilde{f}_n$.\\
 Let $e_n$ be a simple closed geodesic and let $\Tilde{e}_n$ be a lift of $e_n$, let $\Tilde{a}_n$ be the arc $\Tilde{e}_n \cap B_n$, assume that $\Tilde{a}_n$ connects $\Tilde{l}_n$ to $\Tilde{f}_n$, then we have the inequality (see Figure \ref{schema-proof-jm-cyril-fig}),
 $$i (\Tilde{a}_n,\Tilde{B}^{+}_n) \geq i(\Tilde{\kappa}'_{n},\Tilde{B}^{+}_{n}) - (\sharp \left\{\Tilde{a}_n \cap \left<g_{n} \right>\tilde{\kappa}'_{n} \right\} + 1)(i(l_n,B^{+}_{n}) + i(f_n,B^{+}_{n})).$$
 Where $\sharp X$ denotes the cardinal of the set $X$.\\
 \begin{figure}
    \centering
    \includegraphics[scale=0.7]{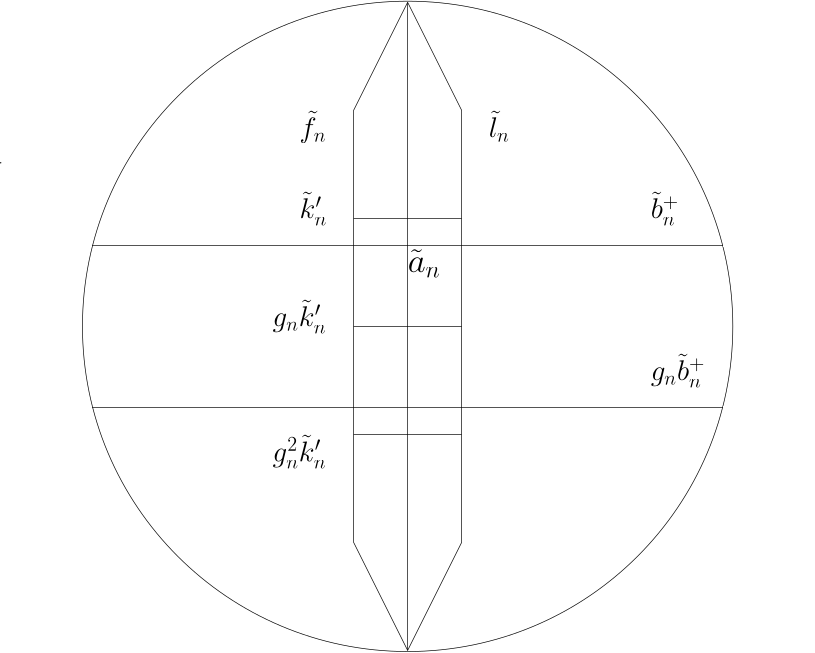}
\caption{Each leaf $\tilde{b}^{+}_{n}$ of $\tilde{B}^{+}_{n}$ that intersects $\tilde{a}_{n}$ but does not intersect $\tilde{k}_{n}$ must intersect the broken geodesics $\tilde{f}_{n}$ and $\tilde{l}_{n}$. This occurs at most $\iota(\tilde{a}_{n}, \tilde{k}'_n) + 1$ times, where $\iota(\tilde{a}_{n}, \tilde{k}'_n)$ denotes the intersection number of $\tilde{a}_{n}$ with the lifts of $\tilde{k}'_n$.}
    \label{schema-proof-jm-cyril-fig}
\end{figure}

Let's show that up to extracting a subsequence, the homotopy classes of $c_n$ don't depend on $n$. Indeed, otherwise we can find a simple closed curve $e$ such that $i(e,c_n) \to \infty$. Since $e$ belongs to a fixed homotopy class we get that $\sharp \left\{a_n \cap \left<g_{n} \right>\tilde{\kappa}'_{n} \right\}$ is bounded (because $\kappa'_{n}$ belongs to $\alpha_n$ which has a free homotopy class independent on $n$). Then the inequality $ \liminf i(e,B^{+}_{n}) \geq \liminf i(e,c_n)\pi$ holds. It follows that $i(e,B^{+}_{n}) \to \infty$, which is absurd, because we assumed that $B^{+}_{n}$ converge, therefore $\liminf i(e,B^{+}_{n})$ must be bounded.\\
In conclusion we can assume that $c_n$ is in a fixed homotopy class.\\
Then, for any simple closed geodesic $e$ that is not equal to $c$ the inequality $\liminf i(e,B^{+}_n) \geq i(e,c)\pi$ holds.
It follows that $c$ is a leaf of $B_{\infty}$ with weight bigger than or equal to $\pi$. The lemma then holds.
\end{proof}
We end the section by giving a proof of Lemma \ref{closing}.
\begin{proof}
 Since by hypothesis of Lemma \ref{closing} the lamination $B_{\infty}$ has no closed leaf of weight bigger than or equal to $\pi$, we deduce that by Lemma \ref{subclos} that $(h^{\pm}_{n})_{n \in \mathbb{N}}$, the sequences of induced metrics on $\partial^{\pm}C(Q_n)$, belong to a compact subset of $\mathcal{T}(S)$. The conclusion follows from Corollary \ref{comsull}.    
\end{proof}

\section{Approximation by first and third fundamental forms}
 In this section we give a proof of Theorem \ref{Main}. Let $Q \in \mathcal{QF}(S)$ be a quasi-Fuchsian manifold, recall that $Q \setminus C(Q)$ is foliated by constant Gaussian curvature surfaces (see Section 2.3). The proof  is based on using Theorem \ref{mix} to realize abstract metrics as first and third fundamental forms of the $k$-surfaces in a sequence of quasi-Fuchsian manifolds. We choose these abstract metrics in a way that they converge to the couple $(h,L) \in \mathcal{ML}(S) \times \mathcal{T}(S)$ that we want to realize. Lemma \ref{closing} will ensure the convergence (up to extracting a sequence) of the sequence of quasi-Fuchsian manifolds. Finally, we conclude by Theorem \ref{extention} that they converge to what we want to realize.\\
In the next lemma we construct the abstract metrics that we will use for approximation.

\begin{lemma}\label{approxi}  Let $(k_{n})_{n \in \mathbb{N}}$ be a sequence of decreasing real numbers that converge to $-1$. Let $h$ be a hyperbolic metric on $S$, and let $\mu$ be a measured lamination on $S$. Then for any $n$ we can find a metric $h_{k_n}$ on $S$ that has constant Gaussian curvature equal to $k_{n}$ such that the sequence $(h_{k_{n}})_{n \in \mathbb{N}}$ converge in the length spectrum to $h$. Also for any $n$ we can find a metric $g_{k_n}$ on $S$ that has constant Gaussian curvature equal to $\frac{k_{n}}{1+k_{n}}$ such that the sequence $(g_{k_{n}})_{n \in \mathbb{N}}$ converges in the length spectrum to $\mu$.
\end{lemma}

By Theorem \ref{mix} we can find a sequence of quasi-Fuchsian manifolds $(Q_{n})_{n \in \mathbb{N}} \subset \mathcal{QF}(S)$, such that the induced metric $I_{k_n}$ on $S^{-}_{k_n}$ is isotopic to $h_{k_n}$ and the third fundamental form $\III_{k_n}$ on $S^{+}_{k_n}$ is isotopic to $g_{k_n}$. Denote the induced metric on $\partial^{-}C(Q_{n})$ by $h^{-}_{n}$ and the bending lamination of $\partial^{+}C(Q_{n})$ by $B^{+}_{n}$. We know that $I_{k_{n}}$ and $\III_{k_{n}}$ converge, but we don't have any information about the convergence of $B^{+}_{n}$ and $h^{-}_{n}$ because for the moment it is unclear whether the sequence $(Q_{n})_{n \in \mathbb{N}}$ converges or not. For this purpose, we need to find upper bounds on the length spectrum of $B^{+}_{n}$ and $h^{-}_{n}$.\\
Recall that $\mathcal{S}$ is the set of free homotopy classes of simple closed curves of $S$ not homotopic to a point.
\begin{lemma}\label{upper}
 Let $\alpha \in \mathcal{S}$, then:
 \begin{itemize}
 \item $\ell_{I_{k_n}}(\alpha) \geq l_{h^{-}_{n}}(\alpha) $.
     \item  $\ell_{\III_{k_{n}}}(\alpha) \geq i(\alpha,B^{+}_{n})$.
     
 \end{itemize}
\end{lemma}
\begin{proof}

By \cite[Lemma II.1.3.4]{canary2006fundamentals},  the projection $r:\Tilde{S}^{-}_{k_n} \to \partial^{-}\Tilde{C}(Q_n)$ is a 1-Lipschitz equivariant map, where $\Tilde{S}^{-}_{k_n}$ and $\partial^{-}\Tilde{C}(Q_n)$ are respectively the lifts of $S^{-}_{k_n}$ and $\partial C^{-}(Q_n)$. Then the first point follows.\\
The second point comes from the existence of an equivarient 1-lipschitz map  from $(\Tilde{S}^{*})^{+}_{k_{n}}$ into $\partial^{+} \Tilde{C}^{*}(Q_{n})$ \cite{AIF_2014__64_2_457_0}, where $(\Tilde{S}^{*})^{+}_{k_{n}}$ and  $\partial^{+} \Tilde{C}^{*}(Q_{n})$ are the dual  of $\Tilde{S}^{+}_{k_n}$ and $\partial^{+}\Tilde{C}(Q_n)$ respectively in de Sitter space. For more details see \cite[Lemma 6.5]{10.2140/gt.2013.17.157}.
\end{proof}
By Lemma \ref{upper}, we deduce that $\forall \alpha \in \mathcal{S}$, the sequence $(l_{h^{-}_{n}}(\alpha))_{n \in \mathbb{N}}$ is bounded, this implies that $(h^{-}_{n})_{n \in \mathbb{N}}$ belongs to a compact subset of $\mathcal{T}(S)$. Also by the same lemma we deduce that the sequence $(i(\alpha,B^{+}_{n}))_{n \in \mathbb{N}}$ is  bounded, this implies that $(B_{n})_{n \in \mathbb{N}}$ converges in the length spectrum (up to extracting a subsequence) to a measured lamination $B_{\infty}$. The next lemma gives us control on the weight of closed leaves of $B_{\infty}$.\\
\begin{lemma}\label{weight}
  Let $\mu$ and $\lambda$ be two measured laminations such that:\begin{itemize}
      \item  $\lambda$ is a discrete lamination in which the weight of every closed leaf is strictly smaller than $\pi$.
      \item $\forall \alpha \in \mathcal{S}$, $i(\alpha,L) \leq i(\alpha,\lambda)$. 
  \end{itemize}
  Then any closed leaf of $\mu$ has weight strictly smaller than $\pi$.
\end{lemma}
\begin{proof}
 Let $\mu'$ be the sublamination of $\mu$ that consists of closed simple geodesics, then the following inequality holds: $$\forall \alpha \in \mathcal{S}, \ i (\alpha,L') \leq i(\alpha,L) \leq i(\alpha,\lambda).$$
 Let $\gamma_{1},..,\gamma_{k}$ be the leaves of $\lambda$, let  $\gamma_{k+1},...,\gamma_{k'}$ be simple closed curves such that $\left\{\gamma_{1},...,\gamma_{k},\gamma_{k+1},...,\gamma_{k'}   \right\}$ forms pants decomposition (it may be that $\lambda$ is already maximal, then we don't need to add leaves).\\
 Note that for any $j$, $i(\gamma_{j},\lambda) = 0$, this implies that $i(\gamma_{j},L') = 0$, so for any $j$, $\gamma_{j}$ does not transversely intersect the leaves of $\mu'$. Let $\beta$ be a leaf of $\mu'$, if $\beta \notin \left\{\gamma_{1},...,\gamma_{k'}   \right\}$ then $\beta$ intersects some $\gamma_{j}$, what is absurd since $\gamma_{j}$ does not transversely intersect the leaves of $\mu'$, we deduce that the leaves of $\mu'$ belong to the set $\left\{\gamma_{1},...,\gamma_{k'}   \right\}$. Let $\sigma_{j}$ be  a simple closed curve dual to $\gamma_{j}$, that is $\sigma_{j}$ is a simple closed curve such that $i(\gamma_{j},\sigma_{j}) = 1$ or  $2$ (depending on the position of $\gamma_{j}$) and $i(\gamma_{r},\sigma_{j}) = 0 $ for any $r \neq j$. Note that the weight of $\gamma_{j}$ in $\lambda$ (respectively $\mu'$) is equal to $\frac{i(\sigma_{j},\lambda)}{i(\sigma_{j},\gamma_{j})}$ (respectively  $\frac{i(\sigma_{j},L')}{i(\sigma_{j},\gamma_{j})}$), also note that if $\gamma_{j}$ does not belong to the lamination then its weight is equal to $0$.\\
 Since $\frac{i(\sigma_{j},L')}{i(\sigma_{j},\gamma_{j})}$ $\leq$ $\frac{i(\sigma_{j},\lambda)}{i(\sigma_{j},\gamma_{j})}$, we deduce that any closed leaf of $\mu$ is a closed leaf of $\lambda$ that has a bigger weight. Then the proof of the lemma follows.
\end{proof}
Now we can give a proof for Theorem \ref{Main}.

\begin{proof}
Let's start proving the first point of the theorem.
 Let $(h,\mu) \in \mathcal{T}(S) \times \mathcal{ML}_{\pi}(S))$. We start by the case when $\mu$  is discrete.\\ 
 Let $k_n$ be a sequence of decreasing numbers converging to $-1$, by Lemma \ref{approxi} there exists a sequence of metrics $h_{k_n}$ on $S$ with a Gaussian curvature equal to $k_n$, that converge to $h$ in the length spectrum. Also by Lemma \ref{approxi} there exists metrics $g_{k_n}$ on $S$ with Gaussian curvature equal to $\frac{k_n}{1+k_n}$, that converge to $\mu$ in the length spectrum.\\  
By Theorem \ref{mix} we can find a sequence of quasi-Fuchsian manifolds $(Q_{n})_{n \in \mathbb{N}} \subset \mathcal{QF}(S)$, such that the induced third fundamental form on $S^{+}_{k_{n}}$ is isotopic to $g_{k_{n}}$, and the induced first fundamental form on $S^{-}_{k_{n}}$ is isotopic to $h_{k_{n}}$.\\
We denote the induced metric on $\partial^{-}C(Q_n)$ by $h^{-}_n$, and we denote the bending lamination on $\partial^{+}C(Q_n)$ by $B^{+}_{n}$.
\begin{sublemma}\label{helpII}
  The following assertions are true:
  \begin{itemize}
      \item The sequence $(h^{-}_{n})_{n \in \mathbb{N}}$ belongs to a compact subset of $\mathcal{T}(S)$.
      \item There exists $B_{\infty} \in \mathcal{ML}_{\pi}(S)$, such that there is a subsequence of $(B^{+}_{n})_{n \in \mathbb{N}}$ that converge to $B_{\infty}$ in the length spectrum.
  \end{itemize}
\end{sublemma}
\begin{proof}
  By Lemma \ref{upper}, for any $\alpha \in \mathcal{S}$, $\ell_{I_{k_n}}(\alpha) \geq \ell_{h^{-}_{n}}(\alpha)$, which implies that $(h^{-}_{n})_{n \in \mathbb{N}}$ belongs to a compact subset of $\mathcal{S}$ (recall that $(I_{k_n})_{n \in \mathbb{N}}$ converge in the length spectrum).\\
Also, Lemma \ref{upper} implies that for any $\alpha \in \mathcal{S}$, $\ell_{\III_{k_n}}(\alpha) \geq i(\alpha,B^{+}_{n})$. Since  $(g_{k_n})_{n \in \mathbb{N}}$ converge to $\mu$ in the length spectrum, we deduce that $(B^{+}_{n})_{n \in \mathbb{N}}$ belongs to compact subset of $\mathcal{ML}(S)$. Then we can say that there exists $B_{\infty} \in \mathcal{ML}(S)$ such that $(B^{+}_{n})_{n \in \mathbb{N}}$ converge to $B_{\infty}$ in the length spectrum(up to extracting a subsequence). Moreover, we have that for any $\alpha \in \mathcal{S}$, $i(\alpha,L) \geq i(\alpha,B_{\infty})$. The fact that $\mu \in \mathcal{ML}_{\pi}(S)$ and that $\mu$ is discrete, allows us to use Lemma \ref{weight} to deduce that $B_{\infty} \in \mathcal{ML}_{\pi}(S)$. \\  
\end{proof}
Sublemma \ref{helpII} ensures that we are in the hypothesis of Lemma \ref{closing}. Then up to extracting a subsequence the quasi-Fuchsian manifolds $Q_{n}$ converge to some quasi-Fuchian manifold $Q_{\infty}$.
To conclude, we will show
\begin{sublemma}\label{almostconclusion}
 The induced hyperbolic metric on $\partial^{-}C(Q_{\infty})$ is $h$, and the bending lamination of $\partial^{+}C(Q_{\infty})$ is $\mu$.
\end{sublemma}
\begin{proof}
 By Theorem \ref{extention} the induced metric on $\partial^{-}C(Q_{\infty})$ is the limit of the metrics $I_{k_n}$, so it is equal to $h$.\\
 Also by Theorem \ref{extention} the bending lamination of $\partial^{+}C(Q_{\infty})$ is the limit of the metrics $g_{k_n}$, then it is $\mu$.
\end{proof}
The following sublemma follows :

\begin{sublemma}\label{helpIII}
 The first point of Theorem \ref{Main} is true when $\mu$ is discrete.   
\end{sublemma} 
 To conclude, let $\mu$ be any lamination in $\mathcal{ML}_{\pi}(S)$, recall that the discrete geodesic measured laminations are dense in the set of measured laminations. Let $(\mu_{n})_{n \in \mathbb{N}} \subset \mathcal{ML}_{\pi}(S)$ be a sequence of discrete laminations that converge to $\mu$. Then by Sublemma \ref{helpIII} there is a sequence of quasi-Fuchsian manifolds $(Q_{n})_{n \in \mathbb{N}} \subset \mathcal{QF}(S)$, such that the induced metric on $\partial^{-}C(Q_n)$ is $h$ and the bending lamination of $\partial^{+}C(Q_n)$ is $\mu_n$. We assumed that $\mu_n$ converge to $\mu$ which is in $\mathcal{ML}_{\pi}(S)$, then we are in the hypothesis of Lemma \ref{closing}. We deduce that the sequence $(Q_n)_{n \in \mathbb{N}}$ converge up to a subsequence to $Q_{\infty} \in \mathcal{QF}(S)$. By \cite[Theorem 4.6]{pjm/1102620682}, the bending lamination of $\partial^{+}C(Q_{\infty})$ is the limit of $(\mu_{n})_{n \in \mathbb{N}}$, so it is equal to $\mu$. And by \cite[Corollary 4.4]{pjm/1102620682}, the induced metric on $\partial^{-}C(Q_{\infty})$ is $h$. This finishes the proof of the first point of Theorem \ref{Main}.\\

 The proof of the second point of Theorem \ref{Main} is similar to the proof of the first point.\\
 Let $(c,\mu) \in \mathcal{T}(S) \times \mathcal{ML}_{\pi}(S)$. Assume first that $\mu$ is discrete. Let $(k_n)_{n \in \mathbb{N}}$ be a sequence of decreasing numbers that converge to $-1$. Take $g_{k_{n}}$ to be metrics on $S$, of constant Gaussian curvature equal $\frac{k_n}{1+k_n}$, and assume that $(g_{k_{n}})_{n \in \mathbb{N}}$ converge to $\mu$. By Theorem \ref{mix} there exists a sequence of quasi-Fuchsian manifolds $(Q_{n})_{n \in \mathbb{N}} \subset \mathcal{QF}(S)$ such that the third fundamental form $\III_{k_n}$ on $S^{+}_{k_n} \subset Q_{n} $ is isotopic to $g_{k_{n}}$ and the induced conformal metric on $\partial^{-}_{\infty}(Q_{\infty})$ is $c$. Denote the induced metric on $\partial^{-}C(Q_n)$ by $h^{-}_{n}$. By Theorem \ref{sullivan} $(h^-_{n})_{n \in \mathbb{N}}$ is in a compact subset of $\mathcal{T}(S)$.
 \begin{sublemma}
   The following assertions are true:
  \begin{itemize}
      \item The sequence $(h^{-}_{n})_{n \in \mathbb{N}}$ belongs to a compact subset of $\mathcal{T}(S)$.
      \item There exists $B_{\infty} \in \mathcal{ML}_{\pi}$, such that there is a subsequence of $(B^{+}_{n})_{n \in \mathbb{N}}$ that converge to $B_{\infty}$ in the length spectrum.
  \end{itemize}  
 \end{sublemma}
 \begin{proof}
   The first point is a direct consequence of Theorem \ref{sullivan}. The second point was already shown in Sublemma \ref{helpII}. 
 \end{proof}
 We have shown the conditions of Lemma \ref{closing} are satisfied. Then up to extracting a subsequence, $(Q_{n})_{n \in \mathbb{N}}$ converge to $Q_{\infty} \in \mathcal{QF}(S)$.
\begin{sublemma}\label{almostconclusion}
 The conformal structure on $\partial_{\infty}^{-}Q_{\infty}$ is $c$, and the bending lamination of $\partial^{+}C(Q_{\infty})$ is $\mu$.
\end{sublemma}
\begin{proof}
  The first point follows from Theorem \ref{bers}, particularly from the fact that the map $B: \mathcal{QF}(S) \to \mathcal{T}(S) \times \mathcal{T}(S)$ is continuous. The second point comes from the fact that the bending lamination on $ \partial^{+}C(Q) $ varies continuously with $Q$ (see \cite[Theorem 4.6]{pjm/1102620682}). 
\end{proof}
The following sublemma follows,
\begin{sublemma}\label{helpIV}
 Let $(c,\mu) \in \mathcal{T}(S) \times \mathcal{ML}_{\pi}(S)$ such that $\mu$ is discrete. Then, there exists a quasi-Fuchsian manifold such that the conformal structure on $\partial^{-}_{\infty}Q$ is $c$, and the the bending lamination of $\partial^{+}C(Q)$ is $\mu$.   
\end{sublemma}
To finish the proof, let $\mu \in \mathcal{ML}_{\pi}(S)$, take a sequence $(\mu_{n})_{n \in \mathbb{N}} \subset \mathcal{ML}_{\pi}(S)$ of discrete laminations that converge to $\mu$. By Sublemma \ref{helpIV} we can find a sequence of quasi-Fuchsian manifolds $(Q_{n})_{n \in \mathbb{N}}$ such that each $Q_n$ induces $c$ as a conformal structure on $\partial_{\infty}^{-}Q_n$, and induces $\mu_n$ as the bending lamination on $\partial^{+}C(Q_n)$. By Lemma \ref{closing} $(Q_{n})_{n \in \mathbb{N}}$ converge, up to subsequence, to a quasi-Fuchsian manifold $Q_{\infty}$. By \cite[Theorem 4.6]{pjm/1102620682}, the bending lamination of $\partial^{+}C(Q_{\infty})$ is $\mu$, and by Theorem \ref{bers} the conformal structure on $\partial_{\infty}^+Q_{\infty}$ is $c$. 
 
\end{proof}

\section{The proof of Theorem A and Theorem A*}
In this section, we will provide a proof of Theorem A and Theorem A*. Let us first recall the statement of the theorems
\begin{theorem}\label{saveme}
Let $h$ and $h^{*}$ respectively be Riemannian metrics on $S$, and denote their curvatures by $k_{h}$ and $k_{h^*}$ respectively (which are not necessarily constant). We assume that $-1 < k_{h}$ and $k_{h^{*}} < 1$. Moreover, we assume that every contractible closed geodesic with respect to $h^{*}$ has length strictly bigger than $2\pi$. Let $\mu$ be a measured lamination on $S$ such that every closed leaf has weight strictly smaller than $\pi$. Then, there exists a convex hyperbolic metric $g$ on $M=S\times(0,1)$, the interior of $\overline{M}=S\times[0,1]$, such that:
\begin{itemize}
    \item $g$ induces the first fundamental form (respectively third fundamental form) on $S\times \left\{ 0\right\}$ which is isotopic to $h$ (respectively $h^{*}$).
    \item $g$ induces on $S \times \left\{ 1\right\}$ a pleated surface structure in which its bending lamination is $\mu$.
\end{itemize}
Furthermore, the surface $ S \times \left\{ 0\right\}$ is smoothly embedded in $M$.
\end{theorem}
Before proceeding with the proof, we need to state some theorems and lemmas that we will use. We begin with the following theorem, proved by Chen and Schlenker.
\begin{theorem}\cite[Theorem 1.6]{chen2022geometric}\label{usethur}
Let $h$ be a Riemannian metric on $S$ of curvature $k_{h}$, that satisfies $-1 < k_{h}$, and let $h^{*}$ be a Riemannian metric on $S$ of curvature $k_{h^{*}} < 1$ and every contractible closed geodesic with respect to $h^{*}$ has length strictly bigger than $2\pi$. Then there exists a unique convex hyperbolic metric $g$ on $M = S \times (0,1)$ such that:
\begin{itemize}
    \item The induced metric on $S \times \left\{0 \right\}$ is isotopic to $h$.
    \item The induced third fundamental form on $S \times \left\{1 \right\}$ is isotopic to $h^{*}$.
\end{itemize} 
\end{theorem}
The following corollary is a direct consequence of Theorem \ref{usethur} and \cite[Proposition 8.3.2]{thurston2022geometry}.
\begin{corollary}\label{usethur2}
Let $h$ be a Riemannian metric on $S$ of curvature $k_{h}$, that satisfies $-1 < k_{h}$, and let $h^{*}$ be a Riemannian metric on $S$ of curvature $k_{h^{*}} < 1$ such that every contractible closed geodesic with respect to $h^{*}$ has length strictly bigger than $2\pi$. There exists a unique quasi-Fuchsian manifold $Q$ (up to isotopy), and a unique geodesically convex domain $C_{h,h^{*}}(Q)$ with boundary such that:
 \begin{itemize}
     \item The inclusion map $\iota : C_{h,h^{*}}(Q) \to Q$ is homotopically equivalent to the identity.
     \item  The boundary of $C_{h,h^{*}}(Q)$ consists of two components $\partial^{+}C_{h,h^{*}}(Q)$ and $\partial^{-}C_{h,h^{*}}(Q)$. The induced metric on $\partial^{-}C_{h,h^{*}}(Q)$ is isotopic to $h$, and the induced third fundamental form on $\partial^{+}C_{h,h^{*}}(Q)$ is isotopic to $h^{*}$.   
 \end{itemize} 
\end{corollary}
As in the proof of Theorem \ref{Main}, we will use Theorem \ref{usethur} to construct a sequence of quasi-Fuchsian manifolds, and we will show that this sequence converges to a quasi-Fuchsian manifold that satisfies Theorem \ref{saveme}.\\
In \cite{MR2208419}, Schlenker defines two maps, $\mathcal{F}$ and $\mathcal{F}^*$. The map $\mathcal{F}$ associates to a convex hyperbolic metric on a manifold $M$ the induced metric (that is, the first fundamental form) on $\partial M$, the boundary of $M$. The map $\mathcal{F}^*$ associates to a convex hyperbolic metric on $M$ the third fundamental form induced on $\partial M$. In our case, $M = S \times (0,1)$ and $\partial M = S \times \{0,1\}$. The main theorem of \cite{MR2208419} shows that the maps $\mathcal{F}$ and $\mathcal{F}^*$ are homeomorphisms, this is concluded in \cite[Section 8]{MR2208419}. For our purposes, we only use the fact that the maps $\mathcal{F}$ and $\mathcal{F}^*$ are continuous, so we formulate the following two lemmas.
\begin{lemma}\label{labour}
 Let $(g_n)_{n\in\mathbb{N}}$ be a sequence of hyperbolic metrics on $M$ with smooth, strictly convex boundary. Let $(h_n)_{n\in\mathbb{N}}$ be the sequence of the first fundamental forms on the boundary. If $(g_n)_{n\in\mathbb{N}}$ converges to a convex hyperbolic metric $g$ on $M$ and if $(h_{n})_{n \in \mathbb{N}}$ converge to a smooth metric $h$ on $\partial M$ with curvature strictly bigger than $-1$. Then the induced first fundamental form on $\partial M$ by $g$ is $h$.
\end{lemma}

\begin{lemma}\label{labour*}
 Let $(g_n)_{n\in\mathbb{N}}$ be a sequence of hyperbolic metrics on $M$ with smooth, strictly convex boundary. Let $(h^{*}_{n})_{n\in\mathbb{N}}$ be the sequence of the third fundamental forms on the boundary. If $(g_n)_{n\in\mathbb{N}}$ converges to a convex hyperbolic metric $g$ on $M$ and if $(h^{*}_{n})_{n \in \mathbb{N}}$ converges to a smooth metric $h^{*}$ on $\partial M$ with curvature strictly smaller than $-1$ with resspect to which every contractible closed geodesic has length bigger then $2\pi$. Then the induced third fundamental form on $\partial M$ by $g$ is $h^{*}$.
\end{lemma}
   
Also, we will need a similar lemma to Lemma \ref{upper}.

\begin{lemma}\label{upper2}
Let $h, h^{*}, Q,$ and $C_{h,h^{*}}(Q)$ be as given in Corollary \ref{usethur2}. The projection map $r:\partial^{-} C_{h,h^{*}}(Q)\rightarrow \partial^{-} C(Q)$ is a  1-Lipschitz map homotopic to the identity.  
\end{lemma}
\begin{proof}
 Let $\Tilde{\partial}^{-}C_{h,h^{*}}(Q)$ and $\Tilde{\partial}^{-} C(Q)$ be lifts of $\partial^{-} C_{h,h^{*}}(Q)$ and $\partial^{-} C(Q)$ respectively.\\ 
 As it was shown in \cite[Lemma II.1.3.4]{canary2006fundamentals} the projection map $r: \Tilde{\partial}^{-}C_{h,h^*}(Q) \to \Tilde{\partial}^{-} C(Q)$ is a 1-Lipschitz equivariant map. Therefore, the statement follows.  
\end{proof}

To prove Theorem A for the third fundamental form $h^{*}$, we will need to show that, given the Gaussian curvature of the third fundamental form on $\partial M$, the first fundamental form is bilipschitz with respect to the third fundamental form as a function of this curvature. To do this we will use some compactness statements given by Labourie \cite{labourie1989immersions} and we use similar techniques as in \cite{Bonsante_2021}.\\
Let's recall statements that we will use.
\begin{theorem}\cite[Theorem D]{labourie1989immersions}\label{D}
 Let $f_{n}:S \rightarrow \mathbb{H}^3$ be a sequence of immersions of a surface $S$ such that the pullbacks $f_{n}^{*}(g_{\mathbb{H}^3})$ of the hyperbolic metric $g_{\mathbb{H}^3}$ converge smoothly to a metric $h$. If
 \begin{itemize}
     \item $(f_{n})_{n \in \mathbb{N}}$ converge in $C^{0}$ to a map $f$.
     \item There exists $k_{0} > -1$ such that for any $n$ the Gaussian curvature of $f_{n}^{*}(g_{\mathbb{H}^3})$ is bigger than $k_{0}$. 
     \item The integral of the mean curvature is uniformly bounded.
 \end{itemize}
 
 Then a subsequence of $(f_{n})_{n \in \mathbb{N}}$ converges smoothly to an isometric immersion $f$ such that $f^{*}(\mathbb{H}^3)=h$.   
\end{theorem}
Next, we introduce the following definition.
\begin{definition}
Let $f : \mathbb{D} \to \mathbb{H}^3$ be an immersion. We say that $f$ is a convex embedding if it is an embedding and if its image is contained in the boundary of a convex subset. 
\end{definition}
We then also state the following lemma.
\begin{lemma}\cite[Lemma 3.7]{Bonsante_2021}\label{bdms}
 Let $f:\mathbb{D}  \rightarrow \mathbb{H}^3$ be a convex embedding and $R$ be the extrinsic diameter of $f(\mathbb{D})$. Denote by $H$ the mean curvature and by $da$ the area form induced by $f$. Then we have 
 $$ \int _{\mathbb{D}}Hda < \frac{A(R+1)}{sinh(1)}, $$
 where $A(\rho)$ denotes the area of the sphere of radius $\rho$ in the hyperbolic space.
\end{lemma}
Now we prove the following :
\begin{lemma}\label{bi}
Let $(M,g)$ be a convex hyperbolic manifold, and let $h^{*}$ be the third fundamental form induced on $S \times \left\{0 \right\}$ (one component of the boundary of $M$). Then there exist $r$ and $r'$, both greater than $0$ and depending only on $h^{*}$, such that the principal curvatures of $S \times \left\{0 \right\}$ belong to the interval $\left [r,r'  \right ]$.\\
In particular, the first fundamental form of $S \times \left\{0 \right\}$ is bi-Lipschitz to the third fundamental form, where the bi-Lipschitz constant depends only on $h^{*}$.  
\end{lemma}
\begin{proof}
 Recall from Section 1.4 that the Gaussian curvature of the first fundamental form of $S \times \left\{ 0 \right\}$, denoted by $k_{I}$, satisfies $k_{I} = \frac{k_{h^{*}}}{1 - k_{h^{*}}}$, where $k_{h^{*}}$ is the curvature of $h^{*}$. Since $S$ is compact and $k_{h^{*}} < 1$, we can deduce that $k_{I}$ is bounded away from $-1$ and $\infty$ by constants that depend only on $k_{h^{*}}$, therefore it depends only on $h^{*}$.\\
We will now argue by contradiction. Assume that the lemma is false. Then, there exists a sequence of convex hyperbolic manifolds $(M,g_n)_{n \in \mathbb{N}}$ such that the induced third fundamental form on $S \times \left\{ 0 \right\}$ is isotopic to $h^*$, and there is a sequence of points $x_n \in (S \times \left\{ 0 \right\},  g_n\mid_{S \times \left\{ 0\right\}})$ such that the principal curvatures at $x_n$ diverge. (Note that their product is equal to $k_I + 1$, so one goes to zero and the other goes to infinity).\\
Let $(D,\tilde{h}_{n})$ be a covering space  of $(S \times \left\{ 0 \right\},  g_n\mid_{S \times \left\{ 0\right\}})$ and let $\Gamma_{n}$ be a quasi-Fuchsian representation (that is $\Gamma_{n}(\pi_{1}(S))$ is a quasi-Fuchsian group) such that $(M,g_{n})$ is embedded in $\mathbb{H}^{3}/\Gamma_{n}$. Then there is an equivariant isometric embedding $f_{n}: (D,\tilde{h}_{n}) \to \mathbb{H}^{3}$ that descends to $(S \times \left\{ 0 \right\},  g_n\mid_{S \times \left\{ 0\right\}})$.\\
Up to applying isometries of $\mathbb{H}^3$, we can assume that $x_{n}$ lifts to a constant $\tilde{x}_{0} \in \mathbb{H}^{3}$.\\
Since $(f_{n})_{n \in \mathbb{N}}$ are isometric embeddings, they are 1-Lipschitz and hence equicontinuous. Up to normalization, we can assume the existence of $o \in D$ such that $f_n(o) = \tilde{x}_0$ for any $n$. Then by the Ascoli-Arzela theorem, up to extracting a subsequence, $(f_{n})_{n \in \mathbb{N}}$ converges in $C^{0}$ topology to some function $f$.\\
Since $k_{\Tilde{h}_{n}}$, the curvature of $\Tilde{h}_{n}$, belongs to a compact interval (independent of $n$), then, up to extracting a subsequence, $\Tilde{h}_{n}$ converge to a metric $\Tilde{h}$.\\
Note that the Gaussian curvature of $\Tilde{h}_{n}$ is bounded away from $-1$ independently on $n$. Also by Lemma \ref{bdms} there is a neighborhood $U$ of the fixed point $o$ such that the integral of the mean curvatures of $f_{n}(U)$ is bounded. Then by Theorem \ref{D} $f_{n} \mid_{U}$ converge to $f \mid_{U}$ in $C^{\infty}$ topology, and moreover, the metrics $(f_{n} \mid_{U})^{*}(g_{\mathbb{H}^3})$ converge to  $(f \mid_{U})^{*}(g_{\mathbb{H}^3})$. In particular, this contradicts the fact that the principal curvature of $f_{n}(U)$ at $x_{0}$ diverge.
\end{proof}

To prove Theorem \ref{saveme}, we will follow a similar approach to the one used in the proof of Theorem \ref{Main}. Namely, we will construct a sequence of quasi-Fuchsian manifolds $(Q_n)_{n \in \mathbb{N}}$ and show that it converges (up to extracting of a sub-sequence). Proposition \ref{labour} and Theorem \ref{extention} ensure that the limit quasi-Fuchsian manifold satisfies the statement of Theorem \ref{saveme}.\\
Let us start the proof.
\begin{proof}
We start by proving the statement for the first fundamental form.\\
Let $h$ be a Riemannian metric on $S$ that has curvature $k_{h}$ strictly bigger than -1. Let $(k_n)_{n \in \mathbb{N}}$ be a sequence of decreasing real numbers that converge to $-1$. Let $\mu$ be a measured lamination on $S$ in which every closed leaf has weight strictly smaller than $\pi$. Let $h_n^*$ be a Riemannian metric on $S$ with a constant Gaussian curvature $k_{h^{*}_n}$ satisfying $k_{h^{*}_n} = \frac{k_n}{1+k_n}$, and such that the sequence $(h^{*}_{n})_{n \in \mathbb{N}}$ converges to $\mu$ in the length spectrum.\\
 Using Corollary \ref{usethur2}, we construct a sequence of quasi-Fuchsian manifolds $(Q_n)_{n \in \mathbb{N}}$, each containing a geodesically convex subset $C_{h,h_n^*}(Q_n)$ as described in the corollary (Corollary \ref{usethur2}).\\
\begin{sublemma}\label{saveme1}
Assume that $\mu$ is a discrete measured lamination, then the sequence $(Q_n)_{n \in \mathbb{N}}$ converges, up to extracting a subsequence, to a quasi-Fuchsian manifold $Q_{\infty}$. 
\end{sublemma}
\begin{proof}
Using Lemma \ref{upper2}, the hyperbolic metrics on $\partial^{-}C(Q_n)$ belong to a compact subset of $\mathcal{T}(S)$. By the second point of Lemma \ref{upper}, we deduce that the bending laminations on $\partial^{+}C(Q_n)$ converge (up to extraction) to a lamination $B_{\infty}$. By Lemma \ref{weight}, $B_{\infty}$ has no closed leaf of weight greater than or equal to zero.\\
Since the hypotheses of Lemma \ref{closing} are satisfied, we can conclude that the Sublemma holds.
\end{proof}
The next Sublemma is a direct consequence of Theorem \ref{extention} and Proposition \ref{labour}.
\begin{sublemma}\label{saveme3}
If $\mu$ is  discrete, then the bending lamination of $\partial^{+} C(Q_{\infty})$ is equal to $\mu$. Moreover, there exists a strictly convex surface $S^{-}_{\infty}$ that is smoothly embedded in $E^{-}(Q_{\infty})$, such that $S^{-}_{\infty}$ bounds, together with $\partial^{+}C(Q_{\infty})$, a geodesically convex subset of $Q_{\infty}$ in which the inclusion is homotopic to the identity. And the induced metric on $S^{-}_{\infty}$ is isotopic to $h$.
\end{sublemma}
\begin{proof}
 The fact that bending lamination of $\partial^{+} C(Q_{\infty})$ is equal to $\mu$ follows directly from Theorem \ref{extention}.\\
 The second point follows from Proposition \ref{labour}. Specifically, it follows from the fact that each geodesically convex hyperbolic manifold $(M,g)$ is embedded in a unique quasi-Fuchsian manifold, as shown in \cite[Proposition 8.3.2]{thurston2022geometry}.  
\end{proof}
We conclude that :
\begin{sublemma}\label{saveme3}
When $\mu$ is discrete, Theorem \ref{saveme} holds for the first fundamental form part.
\end{sublemma}
Next, suppose that $\mu$ is not necessarily discrete. Let $(\mu_{n})_{n \in \mathbb{N}}$ be a sequence of discrete elements of $\mathcal{ML}{\pi}(S)$ (recall that $\mathcal{ML}{\pi}(S)$ is the set of measured laminations in which every closed leaf has weight strictly smaller than $\pi$) that converges to $\mu$. By Sublemma \ref{saveme3}, we can find a sequence of geodesically convex hyperbolic manifolds $(M,g_{n})$ such that the induced metric on $S \times  \left\{ 0\right\}$ is isotopic to $h$, and $g_{n}$ induces on $S \times  \left\{ 1\right\}$ a structure of pleated surface of bending lamination $\mu_{n}$. Since such a submanifold is embedded in a unique quasi-Fuchsian manifold, we can apply the same argument given in Sublemma \ref{saveme1} to deduce that the metrics $g_{n}$ converge to a hyperbolic metric $g$ in which the induced metric on $S \times  \left\{ 0\right\}$ is isotopic to $h$. Moreover, by the continuity of the bending lamination \cite[Theorem 4.6]{pjm/1102620682}, the metric $g$ induces on $S \times \left\{ 0\right\}$ a structure of pleated surface of bending lamination $\mu$.\\

For the third fundamental form, the statement follows in a similar way. First, we assume that $\mu$ is discrete. Let $(M,g_{n})_{n \in \mathbb{N}}$ be a sequence of convex hyperbolic manifolds in which the induced third fundamental form on $S \times \left\{0\right\}$ is isotopic to $h^{*}$ and the third fundamental form on $S \times \left\{1\right\}$ is isotopic to $h^{*}_{n}$, where $h^{*}_{n}$ has constant Gaussian curvature equal to $\frac{k_{n}}{1+k_{n}}$ and converges to $\mu$ in the length spectrum.\\
By Lemma \ref{bi}, there is a positive constant that depends only on the curvature of $h^{*}$ such that $h^{-}_{n}$, the first fundamental form on $S \times  \left\{0 \right\}$ induced by $g_n$, is bi-Lipschitz to $h^{*}$.\\
Let $Q_{n}$ be the unique quasi-Fuchsian manifold in which $(M,g_n)$ is embedded. By Lemma \ref{upper2}, we deduce that the sequence of induced metrics on $\partial^{-}C(Q_n)$ converges, up to extracting a subsequence. Moreover, by Lemma \ref{upper} and Lemma \ref{weight}, the bending lamination $B^{+}_{n}$ on $\partial^{+}C(Q_n)$ converges, up to extracting a subsequence, to a lamination $B_{\infty}$ without any closed leaf of weight bigger than $\pi$.\\
We deduce, by Lemma \ref{closing}, that the metrics $g_n$ converge, up to extracting a subsequence, to a convex hyperbolic metric $g$. Then, by Lemma \ref{labour*}, the third fundamental form on $S \times \left\{0\right\}$ is isotopic to $h^{*}$, and by Lemma \ref{approxi}, the bending lamination on $S \times \{1\}$ is $\mu$.\\
If $\mu$ is not discrete, we approximate it by a sequence of discrete laminations $(\mu_n)_{n \in \mathbb{N}}$. We can then construct a sequence of convex hyperbolic metrics $(g_n)_{n \in \mathbb{N}}$ on $M$ such that the third fundamental form on $S \times \left\{0\right\}$ is isotopic to $h^{*}$ and the bending lamination on $S \times \left\{1\right\}$ is $\mu_n$. Applying Lemma \ref{bi}, Lemma \ref{upper2}, and Lemma \ref{closing}, we find that the metrics $g_n$ converge,up to extracting a subsequence, to a convex hyperbolic metric $g$, and so $(M,g)$ is the desired hyperbolic convex hyperbolic manifold.
\end{proof}

\section{Uniqueness near the Fuchsian locus}

It was shown by Bonahon in \cite{bonahon2005kleinian} that if we assume that two measured laminations $\mu^{+}$ and $\mu^{-}$ are small enough (for the meaning of small, we refer to see Bonahon's paper \cite{bonahon2005kleinian}, or see the statement of Theorem \ref{mixte unicity near the fuchsian locus} of this section), then there exists a unique quasi-Fuchsian manifold $Q$ that has $\mu^{+}$ (resp $\mu^{-}$) as the bending lamination on $\partial^{+}C(Q)$ (resp on $\partial^{-}C(Q)$).\\ 
In this section, we will establish the uniqueness in the first part of the statement of Theorem \ref{Main}, close to the Fuchsian locus (we will provide a definition for "close to the Fuchsian locus" later in this section). This uniqueness holds when $h$ is considered as the hyperbolic metric on the boundary of the convex core, given that we assume the lamination $\mu$ to be sufficiently small. The precise notion of small enough will be established in the main theorem of this section's statement.\\
The main theorem of this section is:
\begin{theorem}\label{mixte unicity near the fuchsian locus}
For any $\mu \in \mathcal{ML}(S)$ and $h \in \mathcal{T}(S)$, there exists an $\delta_{h,\mu} > 0$ such that for any $0 < t < \delta_{h,\mu}$, there exists a unique quasi-Fuchsian manifold $Q$, such that the hyperbolic metric on $\partial^- C(Q)$ is $h$ and the bending lamination on $\partial^+ C(Q)$ is $t\mu$.
\end{theorem}
We recall that the Teichm\"{u}ller space $\mathcal{T}(S)$ defined in Section 1.1, and the space of quasi-Fuchsian manifolds $\mathcal{QF}(S)$ defined in Section 2, possess the structure of smooth (actually analytic) manifolds.\\
Before beginning the proof we will introduce some preliminaries.
\subsection{The length function}
Let's assume that $\mu$ is a discrete measured lamination, which is defined as a set of homotopy classes of disjoint simple closed curves not homotopic to a point, with each homotopy class being associated with a positive scalar referred to as its weight. Recall that discrete measured laminations are dense in the set of measured laminations (for precise definitions and references, see Section 1.2).\\
Let's denote by $c_{i}$ the representatives of the simple closed curves for the homotopy classes that form $\mu$, and let's denote by $a_{i}$ the associated weight for the homotopy class of $c_{i}$. Let $h$ be a hyperbolic metric. Then we define : 
$$\ell_{h}(L) = \sum_{i}a_{i}\ell_{h}(c_{i}).$$
Here, $\ell_{h}(c_{i})$ denotes the length of the unique simple closed geodesic in the hyperbolic surface $(S,h)$ that is homotopic to $c_{i}$.\\
It was demonstrated in \cite{bonahon1986bouts} that $\ell_{h}(L)$ possesses a unique continuous extension when the measured lamination $\mu$ is not discrete. This extension defines a function :
$$\ell : \mathcal{ML}(S) \times \mathcal{T}(S) \to \mathbb{R}_{+},$$
referred to as the length function.\\
Fixing a measured lamination $\mu$, this defines the map : 
\begin{align*}
    \ell_{\mu}: \mathcal{T}(S) &\to \mathbb{R}\\
     h &\mapsto \ell(\mu, h).
\end{align*}

The following theorem is given in \cite[Theorem 1.2]{kerckhoff1992lines}.
\begin{theorem}\label{crit}
    Let $\mu, \mu' \in \mathcal{ML}(S)$ be two laminations that fill $S$. Then $\ell_{\mu} + \ell_{\mu'}$ has a unique critical point on $\mathcal{T}(S)$, which is necessarily a minimum.
\end{theorem}

Under the hypothesis of Theorem \ref{crit}, we denote the critical point of $\ell_{\mu} + \ell_{\mu'}$ by $\kappa(\mu, \mu')$.

Later, Bonahon established the injectivity of $\kappa$, as the following lemma illustrates
\begin{lemma}\label{injectivecrit}\cite[Lemma 4]{bonahon2005kleinian}
 Let $\mu$, $\mu'$, and $\mu''$ be measured laminations in $\mathcal{ML}(S)$. Then, $$\kappa(\mu, \mu') = \kappa(\mu, \mu''),$$ implies that $\mu' = \mu''$.   
\end{lemma}

\subsection{Tangentiability}
For more details see \cite{bonahon1998variations}. A map $f: U \to V$ between two open subsets of $\mathbb{R}^{n}$ is said to be tangentiable if, for each $x \in U$ and $v \in \mathbb{R}^{n}$, the limit :
$$\frac{\partial}{\partial t}\bigg|_{t=0^{+}}f(x+tv) = \lim\limits_{t \to 0^{+}}\frac{f(x+tv)-f(x)}{t}.$$
exists, and the convergence is locally unifrom on $v$. In other words, a tangentiable map can be seen as a a one sided differentiable map with directional derivatives everywhere.\\
If $f$ is tangentiable, we can define the tangent map of $f$ at $x$ as follows : 
\begin{align*}
T_{x}f:  \mathbb{R}^{n} & \to \mathbb{R}^{n}\\
T_{x}f(v) & \mapsto \frac{\partial}{\partial t}\bigg|_{t = 0^{+}}f(x + tv).
\end{align*}
Note that when $f$ is tangentiable, $T_{x}f$ is both continuous and homogeneous (see\cite[Section 1]{bonahon1998variations}), in other words, for any $\lambda > 0$, we have $T_{x}f(\lambda v) = \lambda T_{x}f(v)$.\\
A tangentiable manifold is locally modeled on $\mathbb{R}^{n}$ and has transition maps that are tangentiable, the smooth manifolds are an examples of tangentiable manifolds. If $M$ is a tangentiable manifold, we can refer to its tangent space as $T_{x}M$. However, it's important to note that $T_{x}M$ has the structure of a cone rather than a vector space. Notably, the space of measured laminations $\mathcal{ML}(S)$ possesses a natural tangentiable manifold structure \cite{penner1992combinatorics} and \cite{thurston2022geometry}.\\ 
The notion of a tangentiable map extends to tangentiable manifolds using local charts. A homeomorphism between two tangentiable manifolds is called a bitangentiable homeomorphism if it and its inverse are tangentiable, and if the tangent maps are homeomorphisms everywhere. Bonahon in \cite[Lemma 4]{bonahon1998variations} gave a criterion for this.
\begin{lemma}
    Let $f: M \to N$ be a homeomorphism between tangentiable manifolds. If $f$ is tangentiable, and all of its tangent maps are injective, then $f$ is a bitangentiable homeomorphism.
\end{lemma}

\subsection{Grafting and quasi-Fuchsian manifolds}
For details, refer to \cite{Kulkarni1994}, \cite{kamishima1992deformation}, and \cite{dumas2009complex}. A projective complex structure on $S$ consists of a maximal atlas on $S$, in which the charts map open sets of $S$ to $\mathbb{CP}^1$, and the transitions are restrictions of M\"{o}bius maps. We denote the set of complex projective structures on $S$ considered up to isotopy by $\mathcal{P}(S)$.\\
Note that if $Q \in \mathcal{QF}(S)$ is a quasi-Fuchsian manifold, then it induces projective structures on each of $ S \times \left\{ 0\right\}$ and  $S \times \left\{ 1\right\} $. Indeed,  $S \times \left\{ 0\right\}$ (resp  $S \times \left\{ 1\right\}$) is identified with $\Omega^{+}_{\Gamma} / \Gamma$ (resp $\Omega^{-}_{\Gamma} / \Gamma$), here we assume that $Q$ is identified with $\mathbb{H}^3 / \Gamma$, where $\Gamma$ is a quasi-Fuchsian group and $\Omega^{\pm}_{\Gamma}$ are the connected components of the domain of discontinuity of $\Gamma$ (see Section 2 for notations and definitions).\\
A natural homeomorphism, called grafting, exists between $\mathcal{P}(S)$ and $\mathcal{ML}(S) \times \mathcal{T}(S)$. We will elucidate this mapping in the situation where the $\mathbb{CP}^{1}$ structure is induced from a quasi-Fuchsian manifold (we will clarify this in the next lines). For a general definition, we refer to the mentioned references (\cite{Kulkarni1994}, \cite{kamishima1992deformation}, and \cite{dumas2009complex}).\\
The space of quasi-Fuchsian manifolds can be embedded in the space of $\mathbb{CP}^{1}$ structures via the following correspondence. Let $Q \in \mathcal{QF}(S)$ be a quasi-Fuchsian manifold. Note that knowing the $\mathbb{CP}^{1}$ structure on $S \times \left\{ 0\right\}$ (resp $S \times \left\{ 1\right\}$ ), then we know the quasi-Fuchsian structure on $Q$ (because we will know the quasi-Fuchsian representation from $\pi_{1}(S)$ to $PSL_{2}(\mathbb{C})$). Then this defines an embedding of the space of quasi-Fuchsian manifolds $\mathcal{QF}(S)$ into the space $\mathcal{P}(S)$ of $\mathbb{CP}^{1}$ structures (the embedding depends if we see the space of $\mathbb{CP}^{1}$ structure on $S \times \left\{ 1\right\}$ or on $S \times \left\{ 0\right\}$).\\
Also, we have an embedding of the space $\mathcal{QF}(S)$ of quasi-Fuchsian manifolds into $\mathcal{ML}(S) \times \mathcal{T}(S)$ by considering the hyperbolic metric and the bending lamination on $\partial^{+}C(Q)$ (or on $\partial^{-}C(Q)$). These two maps coincide with the inverse of the grafting map mentioned above when we embed $\mathcal{QF}(S)$ into $\mathcal{P}(S)$.
$$Pr^{\pm}:\mathcal{QF}(S) \to \mathcal{ML}(S) \times \mathcal{T}(S),$$ 
where the sign $+$(res $-$) is taken when we consider the $\mathbb{CP}^1$ structure on $S \times \left\{ 1\right\}$ (resp on $S \times \left\{ 0\right\}$).\\

Bonahon has shown the following Theorem \cite[Theorem 3]{bonahon1998variations}.
\begin{theorem}\label{bit}
  The inverse of the grafting map $Pr :\mathcal{P}(S) \to \mathcal{ML}(S) \times \mathcal{T}(S)$ is a bitangentiable homeomorphism.  
\end{theorem}
We refer to Section 2 to recall the definitions of $\mathcal{F}(S)$ the Fuchsian manifolds and $\mathcal{QF}(S)$ the Quasi-Fuchsian manifolds. We also recall that $Mix_{-1}$ is the map defined from $\mathcal{QF}(S)$ to $\mathcal{ML}_{\pi}(S) \times \mathcal{T}(S)$, which associates to a quasi-Fuchsian manifold $Q$ the hyperbolic metric on $\partial^{-}C(Q)$ and the bending lamination on $\partial^{+}C(Q)$ (see Section 2.4).\\
To initiate the proof of the main statement of this section, we will first establish the existence of a neighborhood of $\mathcal{F}(S)$ denoted as $V$, in which the map $Mix_{-1} \mid _{V}$ is injective. Later we will show that $(Mix_{-1})^{-1}(Mix_{-1}(V)) = V$.\\

We will employ a strategy similar to the one used in \cite{bonahon2005kleinian}, where Bonahon shows the uniqueness of the map $L$ in a neighborhood of the Fuchsian locus (see Section 2.2 for the definition of the map $L$). That is, when the bending laminations are sufficiently small, they uniquely determine the quasi-Fuchsian manifold.\\
Hereafter, let $\mu \in \mathcal{ML}(S)$ be a measured lamination, and $h \in \mathcal{T}(S)$ be a point in the Teichm\"{u}ller space of $S$. We denote by $\mathcal{B}^{+}(\mu)$ the subset of quasi-Fuchsian manifolds $Q$ that possess a bending lamination equal to $t \mu$ on $\partial^{+}\mathcal{C}(Q)$ for some $t \geq 0$. Similarly, $\mathcal{M}^{-}(h)$ is the subset of quasi-Fuchsian manifolds $Q$ that feature $h$ as the induced metric on $\partial^{-}\mathcal{C}(Q).$\\
We denote the dimension of $\mathcal{T}(S)$ by $\theta$.\\
The subset $\mathcal{B}^{+}(\mu)$ is a submanifold with a boundary, as demonstrated by the following lemma.

\begin{lemma}\label{pleatedman}\cite[Lemma 7]{bonahon2005kleinian}
The set $\mathcal{B}^{+}(\mu)$ is a submanifold of $\mathcal{QF}(S)$ with boundary of dimension $\theta + 1$. Furthermore, the boundary $\partial \mathcal{B}^{+}(\mu)$ is $\mathcal{F}(S)$. (Recall that $\theta$ is the dimension of $\mathcal{T}(S)$ and that $\mathcal{F}(S)$ denotes the space of Fuchsian manifolds).
\end{lemma}

We will also show that $\mathcal{M}^{-}(h)$ is a $C^1$ submanifold of $\mathcal{QF}(S)$ (although we do not know if it is of higher regularity than $C^1$). We can deduce this as a corollary from the following theorem.\\
Denote by $gr_{\mathcal{T}}$ the composition of the projection map $\mathcal{ML}(S) \times \mathcal{T}(S) \to \mathcal{T}(S)$ with the inverse of the grafting map $Pr: \mathcal{P}(S) \to \mathcal{ML}(S) \times \mathcal{T}(S).$
\begin{theorem} \cite[Theorem 3]{bonahon1998variations}
The map : 
$$gr_{\mathcal{T}}: \mathcal{P}(S) \to  \mathcal{T}(S)$$
 is $C^1$.
\end{theorem}
As it has been observed in \cite[Corollary 2.5]{dumas2008projective} the map $gr_{\mathcal{T}}$ is a $C^{1}$ submersion.
\begin{corollary}\label{easily}
The map $gr_{\mathcal{T}}: \mathcal{P}(S) \to  \mathcal{T}(S)$ is a $C^{1}$ submersion.  
\end{corollary}
\begin{proof}
 Denote by $Pr:\mathcal{P}(S) \to \mathcal{ML}(S) \times \mathcal{T}(S)$ the inverse of the grafting map. Note that for each $\lambda \in \mathcal{ML}(S)$, the relation $gr_{\mathcal{T}} \circ Pr^{-1}(\lambda,.) = Id$ holds. Then for each $\lambda
  \in \mathcal{ML}(S)$, the map $Pr^{-1}(\lambda,.)$ is a smooth section (smoothness follows from the proof of \cite[Theorem A]{scannell2002grafting}) of $gr_{\mathcal{T}}$, and since the map $Pr$ is a homeomorphism (in particular surjective), the sections fill up $\mathcal{P}(S)$, so we conclude ($gr_{\mathcal{T}}$ is surjective, and its differentials are surjective at each point).    
\end{proof}
Then a direct consequence of Corollary \ref{easily} is the following.
\begin{corollary}\label{metricman}
 The submanifold $\mathcal{M}^{-}(h)$ is $C^1$ submanifold of $\mathcal{QF}(S)$.   
\end{corollary}
\begin{proof}
 Note that $\mathcal{M}^{-}(h) = gr^{-1}_{\mathcal{T}}(h) \cap \mathcal{QF}(S)$. The statement follows because $gr_{\mathcal{T}}$ is a $C^{1}$ submersion.   
\end{proof}

Bonahon showed a proposition which is crucial in this section. Before stating it, we would like to draw the reader's attention to the fact that when $\mu_{0} = 0 \in \mathcal{ML}(S)$, the tangent space $T_{\mu_{0}}\mathcal{ML}(S)$ in $\mathcal{ML}(S)$ corresponds to the set of measured laminations $\mathcal{ML}(S)$ itself (This was mentioned in \cite{bonahon2005kleinian}, for more details we refer to \cite{ASENS_1997_4_30_2_205_0} ).
\\

As in Section 2.2, we denote by :
$$L: \mathcal{QF}(S) \to \mathcal{ML}(S) \times \mathcal{ML}(S),$$ 
the map that associates to a quasi-Fuchsian manifold $Q$ the bending laminations on $\partial^{+}C(Q)$ and $\partial^{-}C(Q)$. 
\begin{lemma}\label{kerkquasi}\cite[Proposition 6]{bonahon2005kleinian}
 Let $\mu, \nu \in \mathcal{ML}(S)$ be two measured geodesic laminations, and let $t \mapsto q_t$, $t \in [0, T]$, be a differentiable curve in $\mathcal{QF}(S)$, originating from a Fuchsian metric $q_0$, such that the derivative $\frac{\partial}{\partial t}L(q_t)|_{t=0^{+}}$ of the bending measured lamination is equal to $(\mu, \nu) \neq 0$. Then $\mu$ and $\nu$ fill up the surface $S$, and $q_0 \in \mathcal{F}(S) = \mathcal{T}(S)$ (we identify the Fuchsian manifolds with Teichm\"{u}ller space) is equal to the minimum $\kappa(\mu,\nu)$ of the length function $l_{\mu} + l_{\nu} : \mathcal{T}(S) \to \mathbb{R}$ (see Section 6.1 for the definition of $\kappa(\mu,\nu)$). 
\end{lemma}

The following lemma follows :

\begin{lemma}\label{tran}
 Let $(q_t)_{t \in [0,\epsilon )}$ be a differentiable curve in $\mathcal{B}^{+}(\mu)$ which is originating from a Fuchsian point\\ $ q_{0} = h_{0} \in \mathcal{F}(S) = \mathcal{T}(S)$. Let $(h_{t},l_{t})$ be respectively the metric and the bending lamination on $\partial^{-}C(q_{t})$. Recall that $q_t \in \mathcal{B}^{+}(\mu)$, which means that the bending lamination on the other side of the convex core , $\partial^{+}C(q_{t})$, is of the form $\alpha(t) \mu$, where $\alpha(t) \geq 0$. Then there exists $\mu' \in \mathcal{ML}(S)$, that depends only on $\mu$ and $q_{0}$, such that $\frac{d}{dt}(l_{t},h_{t})|_{t = 0}$ is of the form $(k \mu',\dot{h})$ where $k \in \mathbb{R}$ and $\dot{h} \in T_{h_{0}}\mathcal{T}(S)$.     
\end{lemma}
\begin{proof}
 First, consider $(\mu_{t},m_{t})$ to be the metric and the bending lamination on $\partial^{+}C(q_{t})$, respectively. As $q_{t}$ belongs to $\mathcal{B}^{+}(\mu)$, there exists a real function (with values in $\mathbb{R}$) denoted as $\alpha(t)$ such that $(\mu_{t},m_{t}) = (\alpha(t)\mu,m_{t})$. This leads to $\frac{d}{dt}(\alpha(t)\mu,m_{t})|_{t = 0} = (\dot{\alpha}\mu, \dot{m})$.

Using Lemma \ref{crit} and Lemma \ref{injectivecrit}, we establish the existence of a unique $\mu' \in \mathcal{ML}(S)$ such that $q_{0} = h_{0} = \kappa(\mu,\mu')$.

Now, let $(l_{t},h_{t})$ represents the metric and the bending lamination on $\partial^{-}C(q_{t})$, respectively. Applying Lemma \ref{kerkquasi}, we find $\dot{l} = \dot{\alpha} \mu'$. Consequently, we deduce that the differential of the coordinates on $\partial^{-} C(q_{t})$ for any differentiable curve in $\mathcal{B}^{+}(\mu)$ takes the form $(k \mu',\dot{h})$, where $\dot{h} \in T_{h_{0}}\mathcal{T}(S)$ and $k$ is a non-negative scalar.
\end{proof}
Now we can deduce that the intersection between $\mathcal{M}^{-}(h)$ and $\mathcal{B}^{+}(\mu)$ at the Fuchsian locus is transverse.\\

\begin{lemma}\label{trans}
The intersection $\mathcal{M}^{-}(h) \cap \partial \mathcal{B}^{+}(\mu)$ consists of exactly one point, namely $h \in \mathcal{F}(S)$. Moreover the intersection $T_{h}\mathcal{M}^{-}(h) \cap T_{h}\mathcal{B}^{+}(\mu)$ consists of exactly one line.  
\end{lemma}
\begin{proof}
The intersection $\mathcal{M}^{-}(h) \cap \partial \mathcal{B}^{+}(\mu)$ consists of the Fuchsian manifold that has $h$ as the induced metric on the convex core (which is a totally geodesic surface in this case). There is a unique such Fuchsian manifold, namely $h \in \mathcal{F}(S)$. Let $q^{+}_{t}$ be a differentiable curve in $\mathcal{B}^{+}(\mu)$ and let $q^{-}_{t}$ be a differentiable curve in $\mathcal{M}^{-}(h)$.\\

Recall that $Pr^{-}$ is the map that associates to a quasi-Fuchsian manifold $Q$ the metric and the bending lamination on $\partial^{-}C(Q)$. This map is bi-tangentiable into its image, as Theorem \ref{bit} states. As we have seen in Lemma \ref{tran}, 
$$
\left.\frac{d}{dt}Pr^{-}(q^{+}_{t})\right|_{t = 0},
$$
is of the form $(k\mu',\dot{m})$, where $\mu'$ is the measured lamination that depends only on $\mu$ and $q_{0} = h$, as explained in Lemma \ref{tran}, and $\dot{m} \in T_{h}\mathcal{T}(S)$ . Also, since the metrics on $\partial^{-}\CH(q^{-}_{t})$ are constant, we have that:
$$
\left.\frac{d}{dt}Pr^{-}(q^{-}_{t})\right|_{t = 0},
$$
is of the form $(\lambda,0)$, where $\lambda \in \mathcal{ML}(S)$.\\

Therefore, the intersection :
$$
T_{h}\mathcal{M}^{-}(h) \cap T_{h}\mathcal{B}^{+}(\mu),
$$ 
is the line generated by :
$$
(d_{h}Pr^{-})^{-1}(\mu',0).
$$

\end{proof}

Next, we deduce the following:
\begin{corollary}\label{submani}
There exists a neighborhood $U_{h,\mu}$ of $h \in \mathcal{F}(S)$ such that the intersection $U_{h,\mu} \cap \mathcal{B}^{+}(\mu) \cap \mathcal{M}^{-}(h)$ forms a submanifold with boundary of dimension 1 in $\mathcal{QF}(S)$.   
\end{corollary}
\begin{proof}
According to Lemma \ref{trans}, the two submanifolds $\mathcal{M}^{-}(h)$ and $\mathcal{B}^{+}(\mu)$ intersect transversely at the point $h$. Consequently, there exists a neighborhood $U_{h,\mu}$ of $h \in \mathcal{F} \subset \mathcal{QF}(S)$, such that the intersection $U_{h,\mu} \cap \mathcal{B}^{+}(\mu) \cap \mathcal{M}^{-}(h)$ forms a one dimensional submanifold with boundary $U_{h,\mu} \cap \partial \mathcal{B}^{+}(\mu) \cap \mathcal{M}^{-}(h)$.    
\end{proof}

\begin{corollary}\label{curve}
For any $h \in \mathcal{T}(S)$ and $\mu \in \mathcal{ML}(S)$, there exists $\epsilon_{\mu,h} > 0$ and a $C^1$ embedded curve :
\begin{align*}
q_{h,\mu}: [ 0,\epsilon_{\mu,h} ) \to & \mathcal{QF}(S)\\
t \mapsto & q_{h,\mu}(t),
\end{align*}
such that $Mix_{-1}(q_{h,\mu}(t)) = (t\mu,h)$.\\
Recall that $Mix_{-1}$ is the map that associates to a quasi-Fuchsian manifold $Q$ the hyperbolic metric on $\partial^- C(Q)$ and the bending lamination on $\partial^+ C(Q)$. 
\end{corollary}
\begin{proof}
From Corollary \ref{submani}, there exists a neighborhood $U_{h,\mu}$ of the Fuchsian locus such that $U_{h,\mu} \cap \mathcal{B}^{+}(\mu) \cap \mathcal{M}^{-}(h)$ forms a $C^1$ submanifold with a boundary.

By definition, if $Q \in U_{h,\mu} \cap \mathcal{B}^{+}(\mu) \cap \mathcal{M}^{-}(h)$, then $Mix_{-1}(Q) = (t\mu,h)$ for some $t \geq 0$. This establishes a mapping:
\begin{align*}
\pi_{h,\mu}: U_{h,\mu} \cap \mathcal{B}^{+}(\mu) \cap \mathcal{M}^{-}(h) & \to \mathbb{R}^{+}\\
Q & \mapsto t,
\end{align*}
where $Mix_{-1}(Q) = (t\mu,h)$.\\
The differentiability of the map $\pi_{h,\mu}$ is shown in the proof of \cite[Lemma 7]{bonahon2005kleinian}.\\
Let $q(t)$ be a curve in $U_{h,\mu} \cap \mathcal{B}^{+}(\mu) \cap \mathcal{M}^{-}(h)$ originating from the Fuchsian locus.\\
We will show that if $\frac{\partial}{\partial t}\pi_{h,\mu}(q(t)) = 0$, then $\frac{\partial}{\partial t}q(t) = 0$. Indeed, the variation of the hyperbolic metric on $\partial^- C(q(t))$ is 0 (due to the constancy of the metrics). And the variation of the bending laminations on $\partial^- C(q(t))$ is also 0. This latter assertion is based on our assumption that if the variation of the bending lamination on $\partial^+ C(q(t))$ is 0, then the variation on $\partial^- C(q(t))$ is also 0. This follows from Lemma \ref{kerkquasi} and the fact that the curve $q$ originates from the Fuchsian point $h$. Therefore, using the bitangentiablity of $Pr^{-}$, we can deduce that $\frac{\partial}{\partial t}q(t) = 0$.\\
Then, applying the inverse function theorem, we find that 
$\pi_{h,\mu}: U_{h,\mu} \cap \mathcal{B}^{+}(\mu) \cap \mathcal{M}^{-}(h) \to \mathbb{R}^{+}$ 
is a local diffeomorphism near 
$U_{h,\mu} \cap \partial \mathcal{B}^{+}(\mu) \cap \mathcal{M}^{-}(h)$. 
As a result, there exists $\epsilon_{h,\mu} > 0$ such that 
the map $\pi_{h,\mu}^{-1}:(0,\epsilon_{h,\mu}) \to \mathcal{QF}(S)$ 
defines an embedded $C^1$ curve. 
Note that $\pi_{h,\mu}^{-1}((0,\epsilon_{h,\mu})) = \{ (t\mu,h) \mid 0 \leq t < \epsilon_{h,\mu} \}$. 
By construction, this curve satisfies the statement of the corollary.

\end{proof}
For $h \in \mathcal{T}(S)$ and $\mu \in \mathcal{ML}(S)$, the set $ \left\{(t\mu,h), 0 \leq t < \epsilon_{h,\mu} \right\}$ is denoted as $ [0,\epsilon_{h,\mu}) (\mu,h)$.\\
Note that from Corollary \ref{curve}, we can define a map
$$\phi_{h,\mu}: [0,\epsilon_{h,\mu})(\mu,h) \to \mathcal{QF}(S) $$ 
such that $Mix_{-1} \circ \phi_{h,\mu} = \text{Id}$ (the map $Mix_{-1}$ is the one defined in Section 2.4). The two submanifolds $\mathcal{M}^{-}(h)$ and $\mathcal{B}^{+}(\mu)$ intersect transversely at each point in the image of $\phi_{h,\mu}$.\\
Choose $\epsilon_{h,\mu}$ in a way that the ray $\left [0,\epsilon_{h,\mu}\right) (\mu,h)$ with the aforementioned properties is maximal (i.e., satisfies $Mix_{-1} \circ \phi_{h,\mu} = \text{Id}$ and ensures transverse intersection between $\mathcal{M}^{-}(h)$ and $\mathcal{B}^{+}(\mu)$ along the image of $\phi_{h,\mu}$). For such chosen $\epsilon_{h,\mu}$ that maximizes the ray, we denote $\mathcal{R}_{h,\mu} :=  [0,\epsilon_{h,\mu}) (\mu,h)$.\\

Bonahon has shown that when $\mu_{n} \in \mathcal{ML}(S)$ converge to $\mu$, then $\mathcal{B}^{+}(\mu_{n})$ converges to $\mathcal{B}^{+}(\mu)$ with respect to the $C^{\infty}$ uniform topology on compact subsets.

\begin{lemma}\cite[Lemma 12]{bonahon2005kleinian}\label{Cinfty}
Let $(\mu_{n})_{n \in \mathbb{N}} $ be a sequence of measured laminations in $\mathcal{ML}(S)$. As $(\mu_{n})_{n \in \mathbb{N}} $ converge to $\mu$ with respect to the topology of $\mathcal{ML}(S)$, the submanifolds $\mathcal{B}^{+}(\mu_{n})$ converge to $\mathcal{B}^{+}(\mu)$ in the $C^{\infty}$ uniform topology on compact subsets.
\end{lemma}

Recall that $\mathcal{M}^{-}(h)$ is defined as the reciprocal image of a submersion. This leads to the following lemma

\begin{lemma}\label{C1}
Let $(h_{n})_{n \in \mathbb{N}} \subset \mathcal{T}(S)$ be a sequence of hyperbolic metrics. When $(h_{n})_{n \in \mathbb{N}}$ converges to $h$ in the topology of $\mathcal{T}(S)$ (described in Section 1.1), then $\mathcal{M}^{-}(h_n)$ tends to $\mathcal{M}^{-}(h)$ in the $C^{1}$ uniform topology on compact subsets.
\end{lemma}
\begin{proof}
Recall that the map $gr_{\mathcal{T}}: \mathcal{P}(S) \to \mathcal{T}(S)$ represents a surjective $C^1$ submersion, and that $\mathcal{M}^{-}(h) = gr_{\mathcal{T}}^{-1}({h}) \cap \mathcal{QF}(S)$. Thus, when $h_{n}$ converges to $h$, the sets $\mathcal{M}^{-}(h_n) = gr_{\mathcal{T}}^{-1}({h_{n}}) \cap \mathcal{QF}(S)$ also converge in the $C^{1}$ uniform topology on compact subsets to $\mathcal{M}^{-}(h) = gr_{\mathcal{T}}^{-1}({h}) \cap \mathcal{QF}(S)$.      
\end{proof}
Let $\mathcal{U} \subset \mathcal{ML}(S) \times \mathcal{T}(S)$ be the union of all the sets $\mathcal{R}_{h,\mu}$ as $(\mu,h)$ varies over $\mathcal{ML}(S) \times \mathcal{T}(S)$.\\
Let $\phi: \mathcal{U} \to \mathcal{QF}(S) $ be the map that restricts to $\phi_{h,\mu}$ at each $\mathcal{R}_{h,\mu}$.\\
Our aim next is to show that $\mathcal{U}$ is an open set, and that the map $\phi$ is a homeomorphism into its image.\\

\begin{proposition}\label{phi}
 The set $\mathcal{U}$ is an open neighborhood of $ \left\{ 0\right\} \times \mathcal{T}(S)$ in $\mathcal{ML}(S) \times \mathcal{T}(S)$, and $\phi$ is a homeomorphism from $\mathcal{U}$ into an open neighborhood $V$ of $\mathcal{F}(S)$ in $\mathcal{QF}(S)$.    
\end{proposition}
\begin{proof}
The restriction $\phi_{h,\mu}$ of $\phi$ to $\mathcal{R}_{h,\mu}$ was constructed by considering the transverse intersection of the manifolds $\mathcal{M}^{-}(h)$ and $\mathcal{B}^{+}(\mu)$ in a neighborhood of the Fuchsian locus $\mathcal{F}(S)$. By Lemma \ref{C1} and Lemma \ref{Cinfty}, the manifold $\mathcal{M}^{-}(h)$ depends continuously on $h$ for the topology of $C^1$-convergence, and $\mathcal{B}^{+}(\mu)$ depends continuously on $\mu$ for the topology of $C^{\infty}$-convergence (therefore $C^{1}$-convergence). It follows that $\epsilon_{h,\mu}$ (the length of $\mathcal{R}_{h,\mu} =  [0,\epsilon_{h,\mu}) (\mu,h)$) depends continuously on $(\mu,h)$, and that $\phi_{h,\mu}$ depends continuously on $(\mu,h)$. We conclude that $\mathcal{U}$, the union of $\mathcal{R}_{h,\mu}$, is open in $\mathcal{ML}(S) \times \mathcal{T}(S)$ and that the map $\phi$ is continuous.\\
Recall that $Mix_{-1} \circ \phi = Id$, it follows that $\phi$ is injective. It follows by the Theorem of the invariance of domain that $\phi$ is a homeomorphism into its image.
\end{proof}
Recall that $Mix_{-1}$ is the map that associates to a quasi-Fuchsian manifold $Q$ the hyperbolic metric on $\partial^- C(Q)$ and the bending lamination on $\partial^+ C(Q)$.\\
The following theorem, is a consequence of Proposition \ref{phi}.
\begin{corollary}\label{Mix-injective}
There is a neighborhood $O$ of the Fuchsian locus, such that $Mix_{-1} \mid_{O}$ is injective.
\end{corollary}
\begin{proof}
 Denote $O := \phi(\mathcal{U})$. From Proposition \ref{phi} the map : $$\phi: \mathcal{U} \to O,$$
 is a homeomorphism and $Mix_{-1} \circ \phi = Id$.\\
 It follows that $Mix_{-1} \mid O = \phi^{-1}$, then in particular it is injective.
\end{proof}

Next, let's show the following :

\begin{proposition}\label{proprness}
Let $h \in \mathcal{T}(S)$, $\mu \in \mathcal{ML}(S)$, and $\theta_{n}$ be a sequence of real numbers converging to $0$. Let $Q_n$ be a sequence of quasi-Fuchsian manifolds such that  $Mix_{-1}(Q_n) = (\theta_{n}\mu,h)$. Then $(Q_n)_{n \in \mathbb{N}}$ converges to the Fuchsian manifold $h$.  
\end{proposition}
\begin{proof}
Recall that by Lemma \ref{closing} the map $Mix_{-1}$ is proper.\\
We will show that each subsequence of $(Q_n)_{n \in \mathbb{N}}$ has a subsequence that converge to $h$.\\
Let $(Q_{\psi(n)})_{n \in \mathbb{N}}$ be a subsequence of $(Q_n)_{n \in \mathbb{N}}$. By definition $Mix_{-1}(Q_{\psi(n)})$ converge to $(0,h)$. Since $Mix_{-1}$ is proper, $(Q_{\psi(n)})_{n \in \mathbb{N}}$ has a subsequence that converge to some quasi-Fuchsian manifold $Q_{\psi(\infty)}$, by continuity of $Mix_{-1}$ we get that $Mix_{-1}(Q_{\psi(\infty)}) = (0,h)$, this implies that $Q_{\psi(\infty)} = h$.     
\end{proof}

\begin{lemma}\label{itself}
There is $O'$, a neighborhood of the Fuchsian locus in $\mathcal{QF}(S)$, such that $(Mix_{-1})^{-1}(Mix_{-1}(O')) = O'$.\\
\end{lemma}
\begin{proof}
We argue by contradiction. If the lemma were not true, then we would find a sequence $(Q_{n})_{n \in \mathbb{N}}$ that does not converge to some point in $\mathcal{F}(S)$, while $((Mix_{-1})^{-1}(Q_n))_{n \in \mathbb{N}}$ converge to some point of $ \left\{0 \right\} \times \mathcal{T}(S)$. This contradicts Proposition \ref{proprness}.
\end{proof}

Then we prove Theorem \ref{mixte unicity near the fuchsian locus}.

\begin{proof}
Let $V := O' \cap O$. Note that $\text{Mix}_{-1} \mid_V$ is injective (since $V \subset O$; see Lemma \ref{Mix-injective}), and that $\text{Mix}_{-1}^{-1}(\text{Mix}_{-1}(V)) = V$, by Lemma \ref{itself} and because $V \subset O'$. \\ 
Now, for any $(\mu,h) \in \mathcal{ML}(S) \times \mathcal{T}(S)$, it suffices to choose $\epsilon > 0$ such that for all $0 < t < \epsilon$, we have $(\text{Mix}_{-1})^{-1}(t\mu,h) \subset V$.

\end{proof}

\medskip

\bibliographystyle{alpha}
\bibliography{mybibliography.bib}

\end{document}